\documentclass[a4paper,twoside,10pt]{amsart}
\usepackage{amsmath,amsfonts,amssymb,amsthm}
\newtheorem{theorem}{Theorem}[section]

\usepackage{mathrsfs}
\usepackage{color,graphicx}
\usepackage[a4paper]{geometry}
\numberwithin{equation}{section}

\author[G. Nemes]{Gerg\H{o} Nemes}
\address{Central European University, Department of Mathematics and its Applications, H-1051 Budapest, N\'ador utca 9, Hungary}
\email{nemesgery@gmail.com}

\begin{document}

\title[Resurgence of the Anger--Weber function]{The resurgence properties of\\ the large order asymptotics of\\ the Anger--Weber function II}

\begin{abstract} In this paper, we derive a new representation for the Anger--Weber function, employing the reformulation of the method of steepest descents by C. J. Howls (Howls, Proc. R. Soc. Lond. A \textbf{439} (1992) 373--396). As a consequence of this representation, we deduce a number of properties of the large order asymptotic expansion of the Anger--Weber function, including explicit and realistic error bounds, asymptotic approximations for the late coefficients, exponentially improved asymptotic expansions, and the smooth transition of the Stokes discontinuities.
\end{abstract}
\maketitle

\section{Introduction}\label{section1}

In the first part of this series of papers \cite{Nemes2}, we proved new resurgence-type representations for the remainder term of the asymptotic expansion of the Anger--Weber function $\mathbf{A}_{ - \nu } \left( \lambda \nu \right)$ with complex $\nu$ and $\lambda \geq 1$. These resurgence formulas have different forms according to whether $\lambda >1$ or $\lambda=1$. The main goal of this paper is to derive a similar representation for the Anger--Weber function $\mathbf{A}_{\nu} \left( \lambda \nu\right)$ with complex $\nu$ and $\lambda > 0$. Our derivation is based on the reformulation of the method of steepest descents by Howls \cite{Howls}. Using this representation, we obtain a number of properties of the large order asymptotic expansion of the Anger--Weber function, including explicit and realistic error bounds, asymptotics for the late coefficients, exponentially improved asymptotic expansions, and the smooth transition of the Stokes discontinuities.

Our first theorem describes the resurgence properties of the asymptotic expansion of $\mathbf{A}_{\nu} \left(\lambda \nu \right)$ for $\lambda > 0$. The notations follow the ones given in \cite[p. 298]{NIST}. Throughout this paper, empty sums are taken to be zero.

\begin{theorem}\label{thm1} Let $\lambda>0$ be a fixed positive real number, and let $N$ be a non-negative integer. Then we have
\begin{equation}\label{eq11}
\mathbf{A}_\nu  \left( {\lambda \nu } \right) = \frac{1}{\pi}\sum\limits_{n = 0}^{N - 1} {\frac{{\left( {2n} \right)!a_n \left( {\lambda } \right)}}{\nu ^{2n + 1}}}  + R_N \left( {\nu ,\lambda } \right)
\end{equation}
for $\left|\arg \nu\right| < \frac{\pi}{2}$, with
\begin{equation}\label{eq9}
a_n \left( \lambda  \right) = \frac{1}{{\left( {2n} \right)!}}\left[ {\frac{d^{2n}}{dt^{2n}}\left( {\frac{t}{{\lambda \sinh t + t}}} \right)^{2n + 1} } \right]_{t = 0} = \frac{{\left( { - 1} \right)^n }}{{\left( {2n} \right)!}}\int_0^{ + \infty } {t^{2n} e^{ - \pi t} iH_{it}^{\left( 1 \right)} \left( {\lambda it} \right)dt}
\end{equation}
and
\begin{equation}\label{eq10}
R_N \left( {\nu ,\lambda } \right) = \frac{{\left( { - 1} \right)^N }}{{\pi \nu ^{2N + 1} }}\int_0^{ + \infty } {\frac{{t^{2N} e^{ - \pi t} }}{{1 + \left( {t/\nu } \right)^2 }}iH_{it}^{\left( 1 \right)} \left( \lambda it \right)dt} .
\end{equation}
\end{theorem}

In a previous paper \cite{Nemes}, we proved similar representations for the Hankel function $H_\nu^{\left(1\right)}\left(\lambda \nu\right)$ with $\lambda \geq 1$. In particular, for any non-negative integer $N$ and fixed $0 < \beta <\frac{\pi}{2}$, we have
\begin{equation}\label{eq30}
H_\nu ^{\left( 1 \right)} \left( {\nu \sec \beta } \right) = \frac{{e^{i\nu \left( {\tan \beta  - \beta } \right) - \frac{\pi}{4}i} }}{{\left( {\frac{1}{2}\nu \pi \tan \beta } \right)^{\frac{1}{2}} }}\left( {\sum\limits_{m = 0}^{M - 1} {\left( { - 1} \right)^m \frac{{U_m \left( {i\cot \beta } \right)}}{{\nu ^m }}}  + R_M^{\left(H\right)} \left( {\nu ,\beta } \right)} \right)
\end{equation}
and
\begin{equation}\label{eq31}
H_\nu ^{\left( 1 \right)} \left( \nu  \right) =  - \frac{2}{{3\pi }}\sum\limits_{n = 0}^{N - 1} {d_{2n} e^{\frac{{2\left( {2n + 1} \right)\pi i}}{3}} \sin \left( {\frac{{\left( {2n + 1} \right)\pi }}{3}} \right)\frac{{\Gamma \left( {\frac{{2n + 1}}{3}} \right)}}{{\nu ^{\frac{{2n + 1}}{3}} }}}  + R_N^{\left( H \right)} \left( \nu  \right),
\end{equation}
for $-\frac{\pi}{2} < \arg \nu < \frac{3\pi}{2}$, with
\begin{gather}\label{eq32}
\begin{split}
U_m \left( {i\cot \beta } \right) & = \left( { - 1} \right)^m \frac{{\left( {i\cot \beta } \right)^m }}{{2^m m!}}\left[ {\frac{{d^{2m} }}{{dt^{2m} }}\left( {\frac{1}{2}\frac{{t^2 }}{{i\cot \beta \left( {t - \sinh t} \right) + \cosh t - 1}}} \right)^{m + \frac{1}{2}} } \right]_{t = 0}\\
& = \frac{i^m}{2\left( {2\pi \cot \beta } \right)^{\frac{1}{2}}}\int_0^{ + \infty } {t^{m - \frac{1}{2}} e^{ -t \left( {\tan \beta  - \beta } \right)} \left( {1 + e^{ - 2\pi t} } \right)i H_{it}^{\left( 1 \right)} \left( {it\sec \beta } \right)dt}
\end{split}
\end{gather}
and
\begin{equation}\label{eq33}
d_{2n} = \frac{1}{{\left( {2n} \right)!}}\left[ {\frac{{d^{2n} }}{{dt^{2n} }}\left( {\frac{{t^3 }}{{\sinh t - t}}} \right)^{\frac{{2n + 1}}{3}} } \right]_{t = 0}  = \frac{{\left( { - 1} \right)^n }}{{\Gamma \left( {\frac{{2n + 1}}{3}} \right)}}\int_0^{ + \infty } {t^{\frac{{2n - 2}}{3}} e^{ - 2\pi t} i H_{it}^{\left( 1 \right)} \left( {it} \right)dt} .
\end{equation}
The remainder terms $R_M^{\left(H\right)} \left( {\nu ,\beta } \right)$ and $R_N^{\left( H \right)} \left( \nu  \right)$ can be expressed as
\begin{equation}\label{eq34}
R_M^{\left(H\right)} \left( {\nu ,\beta } \right) =  \frac{1}{{2\left( {2\pi \cot \beta } \right)^{\frac{1}{2}} \left(i\nu\right)^M }}\int_0^{ + \infty } {\frac{{t^{M - \frac{1}{2}} e^{ - t\left( {\tan \beta  - \beta } \right)} }}{{1 + it/\nu }}\left( {1 + e^{ - 2\pi t} } \right)i H_{it}^{\left( 1 \right)} \left( {it\sec \beta } \right)dt}
\end{equation}
and
\begin{equation}\label{eq35}
R_N^{\left( H \right)} \left( \nu  \right) = \frac{{\left( { - 1} \right)^N }}{{3\pi \nu ^{\frac{{2N + 1}}{3}} }}\int_0^{ + \infty } {t^{\frac{{2N - 2}}{3}} e^{ - 2\pi t} \left( {\frac{{e^{\frac{{\left( {2N + 1} \right)\pi i}}{3}} }}{{1 + \left( {t/\nu } \right)^{\frac{2}{3}} e^{\frac{{2\pi i}}{3}} }} + \frac{1}{{1 + \left( {t/\nu } \right)^{\frac{2}{3}} }}} \right)H_{it}^{\left( 1 \right)} \left( {it} \right)dt} .
\end{equation}
These representations of the Hankel function will play an essential role in later sections of this paper. It is important to note that for the case $0 < \lambda < 1$, no simple explicit expression is known for the remainder of the asymptotic series of the Hankel function $H_\nu^{\left(1\right)}\left(\lambda \nu\right)$.

If $\mathbf{J}_{\nu}\left(z\right)$ denotes the Anger function, then $\sin \left( {\pi \nu } \right)\mathbf{A}_{\nu} \left( \lambda \nu  \right) =\mathbf{J}_{\nu} \left(\lambda\nu\right)- J_{\nu} \left(\lambda\nu\right)$ (see \cite[p. 296]{NIST}). From these and the continuation formulas for the Bessel and Hankel functions (see \cite[p. 222 and p. 226]{NIST}), we find
\begin{align*}
\sin \left( {\pi \nu } \right)\mathbf{A}_{\nu} \left( {\lambda \nu e^{2\pi im} } \right) & = \mathbf{J}_{\nu} \left( {\lambda \nu e^{2\pi im} } \right) - J_{\nu } \left( {\lambda \nu e^{2\pi im} } \right)
\\ & = \sin \left( {\pi \nu } \right)\mathbf{A}_\nu  \left( {\lambda \nu } \right) + \left( {1 - e^{2\pi im\nu } } \right)J_\nu  \left( {\lambda \nu } \right)
\\ & = \sin \left( {\pi \nu } \right)\mathbf{A}_\nu  \left( {\lambda \nu } \right) - ie^{\pi im\nu } \sin \left( {\pi m\nu } \right)\left( {H_\nu ^{\left( 1 \right)} \left( {\lambda \nu } \right) + H_\nu ^{\left( 2 \right)} \left( {\lambda \nu } \right)} \right)
\end{align*}
for every integer $m$. From this expression and the resurgence formulas \eqref{eq11}, \eqref{eq30} and \eqref{eq31}, we can derive analogous representations in sectors of the form
\[
\left( {2m - \frac{1}{2}} \right)\pi  < \arg \nu  < \left( {2m + \frac{1}{2}} \right)\pi ,\; m \in \mathbb{Z},
\]
as long as $\lambda \geq 1$. Similarly, applying the continuation formulas
\begin{gather}\label{eq50}
\begin{split}
& - \sin \left( {\pi \nu } \right)\mathbf{A}_{-\nu} \left( {\lambda \nu e^{\left( {2m + 1} \right)\pi i} } \right) = \mathbf{J}_{-\nu} \left( {\lambda \nu e^{\left( {2m + 1} \right)\pi i} } \right) - J_{-\nu} \left( {\lambda \nu e^{\left( {2m + 1} \right)\pi i} } \right)\\
& = \sin \left( {\pi \nu } \right)\mathbf{A}_\nu  \left( {\lambda \nu } \right) + J_\nu  \left( {\lambda \nu } \right) - e^{ - \left( {2m + 1} \right)\pi i\nu } J_{ - \nu } \left( {\lambda \nu } \right)\\
& = \sin \left( {\pi \nu } \right) \mathbf{A}_\nu  \left( {\lambda \nu } \right) + ie^{ - \pi i m \nu} \sin \left( {\pi m\nu } \right)H_\nu ^{\left( 1 \right)} \left( {\lambda \nu } \right) + ie^{ - \left( {m + 1} \right)\pi i\nu } \sin \left( {\pi \left( {m + 1} \right)\nu } \right)H_\nu ^{\left( 2 \right)} \left( {\lambda \nu } \right)
\end{split}
\end{gather}
and the representations \eqref{eq11}, \eqref{eq30} and \eqref{eq31}, we can obtain resurgence formulas in any sector of the form
\[
\left( {2m + \frac{1}{2}} \right)\pi  < \arg \nu  < \left( {2m + \frac{3}{2}} \right)\pi ,\; m \in \mathbb{Z},
\]
provided that $\lambda \geq 1$. The lines $\arg \nu  = \left( {2m \pm \frac{1}{2}} \right)\pi$ are the Stokes lines for the function $\mathbf{A}_{\nu} \left(\lambda \nu\right)$.

When $\nu$ is an integer, the limiting values have to be taken in these continuation formulas.

If we neglect the remainder term and extend the sum to $N = \infty$ in Theorem \ref{thm1}, we recover the known asymptotic series of the Anger--Weber function. Some other formulas for the coefficients $a_n\left(\lambda\right)$ can be found in \cite[Appendix A]{Nemes2}.

In the following two theorems, we give exponentially improved asymptotic expansions for the function $\mathbf{A}_{\nu } \left(\lambda \nu\right)$ when $\lambda>1$ and $\lambda=1$, respectively. Since there is no simple resurgence formula for the Hankel function $H_\nu^{\left(1\right)}\left(\lambda \nu\right)$ when $0<\lambda <1$, at least with our method, we can not prove exponentially improved expansions for the function $\mathbf{A}_{\nu } \left(\lambda \nu\right)$ in the range $0<\lambda <1$. We express our expansions in terms of the Terminant function $\widehat T_p\left(w\right)$ whose definition and basic properties are given in Section \ref{section5}. In these theorems, $R_N \left( {\nu ,\lambda } \right)$ is defined by \eqref{eq11} and it is extended to the sector $\left|\arg \nu\right| \leq \frac{3\pi}{2}$ via analytic continuation. In Theorem \ref{thm2}, we employ the substitution $1 < \lambda = \sec \beta$ with a suitable $0 < \beta <\frac{\pi}{2}$. Throughout this paper, we write $\mathcal{O}_{M,\rho }$ to indicate the dependence of the implied constant on the parameters $M$ and $\rho$.

\begin{theorem}\label{thm2} Suppose that $\left|\arg \nu\right| \leq \frac{3\pi}{2}$, $\left|\nu\right|$ is large and $N = \frac{1}{2}\left| \nu  \right|\left( {\tan \beta  - \beta +\pi} \right) + \rho$ is a positive integer with $\rho$ being bounded. Then
\begin{multline*}
R_N \left( {\nu ,\sec \beta } \right) =  i\frac{{e^{i\nu \left( {\tan \beta  - \beta +\pi} \right) - \frac{\pi }{4}i} }}{{\left( {\frac{1}{2}\nu \pi \tan \beta } \right)^{\frac{1}{2}} }}\sum\limits_{m = 0}^{M - 1} {\left( { - 1} \right)^m \frac{{U_m \left( {i\cot \beta } \right)}}{{\nu ^m }}\widehat T_{2N - m + \frac{1}{2}} \left( {i\nu \left( {\tan \beta  - \beta +\pi} \right)} \right)}\\  - i\frac{{e^{ - i\nu \left( {\tan \beta  - \beta +\pi} \right) + \frac{\pi }{4}i} }}{{\left( {\frac{1}{2}\nu \pi \tan \beta } \right)^{\frac{1}{2}} }}\sum\limits_{m = 0}^{M - 1} {\frac{{U_m \left( {i\cot \beta } \right)}}{{\nu ^m }}\widehat T_{2N - m + \frac{1}{2}} \left( { - i\nu \left( {\tan \beta  - \beta +\pi} \right)} \right)}  + R_{N,M} \left( {\nu ,\sec \beta } \right)
\end{multline*}
with $M$ being an arbitrary fixed non-negative integer, and
\[
R_{N,M} \left( {\nu ,\sec \beta } \right) = \mathcal{O}_{M,\rho } \left( {\frac{{e^{ - \left| \nu  \right|\left( {\tan \beta  - \beta +\pi} \right)} }}{{\left( {\frac{1}{2}\left| \nu  \right|\pi \tan \beta } \right)^{\frac{1}{2}} }}\frac{{\left| {U_M \left( {i\cot \beta } \right)} \right|}}{{\left| \nu  \right|^M }}} \right)
\]
for $\left|\arg \nu\right| \leq \frac{\pi}{2}$;
\[
R_{N,M} \left( {\nu ,\sec \beta } \right) = \mathcal{O}_{M,\rho } \left( {\frac{{e^{ \mp \Im \left( \nu  \right)\left( {\tan \beta  - \beta +\pi} \right)} }}{{\left( {\frac{1}{2}\left| \nu  \right|\pi \tan \beta } \right)^{\frac{1}{2}} }}\frac{{\left| {U_M \left( {i\cot \beta } \right)} \right|}}{{\left| \nu  \right|^M }}} \right)
\]
for $\frac{\pi}{2} \leq \pm \arg \nu \leq \frac{3\pi}{2}$.
\end{theorem}

\begin{theorem}\label{thm3} Suppose that $\left|\arg \nu\right| \leq \frac{3\pi}{2}$, $\left|\nu\right|$ is large and $N = \frac{1}{2}\pi \left| \nu  \right| + \rho$ is a positive integer with $\rho$ being bounded. Then
\begin{multline*}
R_N \left( {\nu ,1} \right) = - ie^{\pi i\nu } \frac{2}{{3\pi }}\sum\limits_{m = 0}^{M - 1} {d_{2m} e^{\frac{{2\left( {2m + 1} \right)\pi i}}{3}} \sin \left( {\frac{{\left( {2m + 1} \right)\pi }}{3}} \right)\frac{{\Gamma \left( {\frac{{2m + 1}}{3}} \right)}}{{\nu ^{\frac{{2m + 1}}{3}} }}\widehat T_{2N - \frac{{2m - 2}}{3}} \left( {\pi i\nu } \right)} 
\\ - ie^{ - \pi i\nu } \frac{2}{{3\pi }}\sum\limits_{m = 0}^{M - 1} {d_{2m} \sin \left( {\frac{{\left( {2m + 1} \right)\pi }}{3}} \right)\frac{{\Gamma \left( {\frac{{2m + 1}}{3}} \right)}}{{\nu ^{\frac{{2m + 1}}{3}} }}\widehat T_{2N - \frac{{2m - 2}}{3}} \left( { - \pi i\nu } \right)}  + R_{N,M} \left( {\nu ,1} \right).
\end{multline*}
with $M$ being an arbitrary fixed non-negative integer, and
\[
R_{N,M} \left( {\nu ,1} \right) = \mathcal{O}_{M,\rho} \left( {e^{ - \pi \left| \nu  \right|} \left| {d_{2M} } \right|\frac{{\Gamma \left( {\frac{{2M + 1}}{3}} \right)}}{{\left| \nu  \right|^{\frac{{2M + 1}}{3}} }}} \right) \; \text{ if } \; M \equiv 0,2 \mod 3,
\]
\[
R_{N,M} \left( {\nu ,1} \right) = \mathcal{O}_{M,\rho} \left( {e^{ - \pi \left| \nu  \right|} \left| {d_{2M + 2} } \right|\frac{{\Gamma \left( {\frac{{2M + 3}}{3}} \right)}}{{\left| \nu  \right|^{\frac{{2M + 3}}{3}} }}} \right) \; \text{ if } \; M \equiv 1 \mod 3
\]
for $\left|\arg \nu\right| \leq \frac{\pi}{2}$;
\[
R_{N,M} \left( {\nu ,1} \right) = \mathcal{O}_{M,\rho } \left( {e^{ \mp \pi \Im \left( \nu  \right)} \left| {d_{2M} } \right|\frac{{\Gamma \left( {\frac{{2M + 1}}{3}} \right)}}{{\left| \nu  \right|^{\frac{{2M + 1}}{3}} }}} \right) \; \text{ if } \; M \equiv 0,2 \mod 3,
\]
\[
R_{N,M} \left( {\nu ,1} \right) = \mathcal{O}_{M,\rho } \left( {e^{ \mp \pi \Im \left( \nu  \right)} \left| {d_{2M + 2} } \right|\frac{{\Gamma \left( {\frac{{2M + 3}}{3}} \right)}}{{\left| \nu  \right|^{\frac{{2M + 3}}{3}} }}} \right) \; \text{ if } \; M \equiv 1 \mod 3
\]
for $\frac{\pi}{2} \leq \pm \arg \nu \leq \frac{3\pi}{2}$.
\end{theorem}

The rest of the paper is organized as follows. In Section \ref{section2}, we prove the resurgence formula stated in Theorem \ref{thm1}. In Section \ref{section3}, we give explicit and realistic error bounds for the asymptotic expansions of $\mathbf{A}_{\nu}\left(\lambda \nu \right)$ using the results of Section \ref{section2}. In Section \ref{section4}, asymptotic approximations for $a_n \left( \lambda \right)$ as $n \to +\infty$ are given. In Section \ref{section5}, we prove the exponentially improved expansions presented in Theorems \ref{thm2} and \ref{thm3}, and provide a detailed discussion of the Stokes phenomenon related to the expansions of $\mathbf{A}_{\nu}\left( \lambda \nu \right)$, $\lambda \geq 1$. The paper concludes with a short discussion in Section \ref{section6}.

\section{Proof of the resurgence formula}\label{section2}

Our analysis is based on the integral definition of the Anger--Weber function
\[
\mathbf{A}_{\nu } \left( z \right) = \frac{1}{\pi }\int_0^{ + \infty } {e^{-\nu t - z\sinh t} dt} \quad \left| {\arg z} \right| < \frac{\pi }{2}.
\]
If $z = \lambda \nu$, where $\lambda$ is a positive constant, then
\begin{equation}\label{eq1}
\mathbf{A}_\nu \left( {\lambda \nu} \right) = \frac{1}{\pi }\int_0^{ + \infty } {e^{ - \nu \left( {\lambda\sinh t + t} \right)} dt} \quad \left| {\arg \nu } \right| < \frac{\pi }{2}.
\end{equation}
The saddle points of the integrand are the roots of the equation $\lambda \cosh t = -1$. Hence, the saddle points are given by $t_{\pm}^{\left( k \right)} = \pm \mathrm{sech}^{ - 1} \lambda + \left(2k+1\right)\pi i$ where $k$ is an arbitrary integer. We denote by $\mathscr{C}_{\pm}^{\left( k \right)}\left(\theta\right)$ the portion of the steepest paths that pass through the saddle point $t_{\pm}^{\left( k \right)}$. Here, and subsequently, we write $\theta = \arg \nu$. As for the path of integration $\mathscr{P}\left(\theta \right)$ in \eqref{eq1}, we take that connected component of
\[
\left\{ {t \in \mathbb{C}:\arg \left[ {e^{i\theta } \left( {\lambda \sinh t + t} \right)} \right] = 0} \right\} \cup \left\{0\right\},
\]
which contains the origin. We remark that $\mathscr{P}\left(0\right)$ is the positive real axis.

\begin{figure}[!t]
\def\svgwidth{0.6\columnwidth}
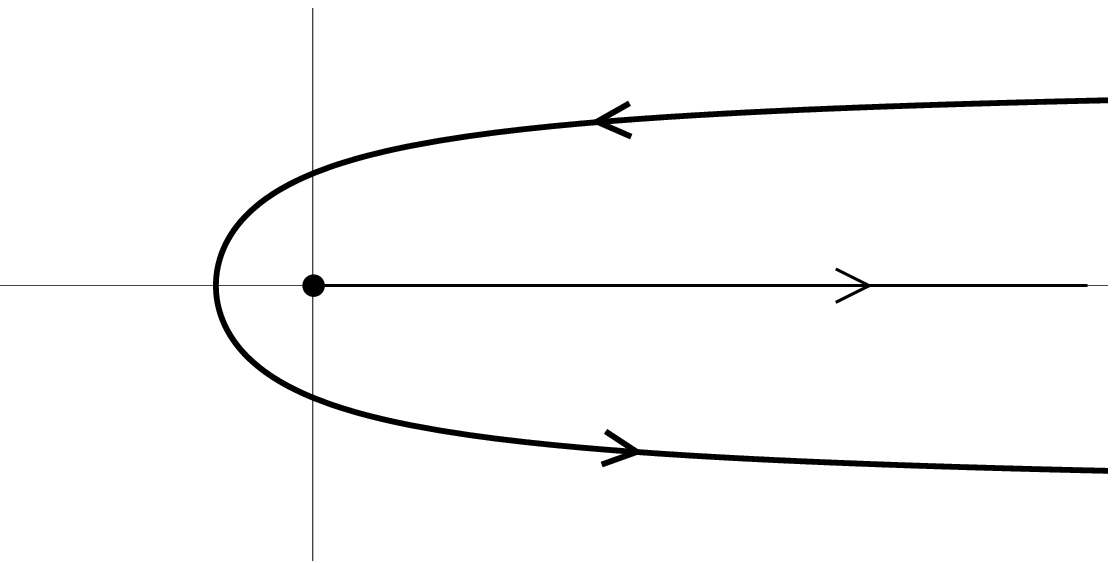
\caption{The contour $\Gamma$ encircling the path $\mathscr{P}\left(0\right)$.}
\label{fig1}
\end{figure}

First, we suppose that $\lambda > 1$ and take $\lambda = \sec \beta$ with a suitable $0<\beta<\frac{\pi}{2}$. With this notation, $t_{\pm}^{\left( k \right)} = \pm i \beta + \left(2k+1\right)\pi i$. For simplicity, we assume that $\theta = 0$. In due course, we shall appeal to an analytic continuation argument to extend our results to complex $\nu$. Let $f\left( {t,\beta } \right) = \sec \beta \sinh t + t$. If
\begin{equation}\label{eq2}
\tau  = f\left( {t,\beta } \right),
\end{equation}
then $\tau$ is real on the curve $\mathscr{P}\left(0\right)$, and, as $t$ travels along this curve from $0$ to $+\infty $, $\tau$ increases from $0$ to $+\infty$. Therefore, corresponding to each positive value of $\tau$, there is a value of $t$, say $t\left(\tau\right)$, satisfying \eqref{eq2} with $t\left(\tau\right)>0$. In terms of $\tau$, we have
\[
\mathbf{A}_\nu  \left( {\nu \sec \beta } \right) = \frac{1}{\pi }\int_0^{ + \infty } {e^{ - \nu \tau } \frac{{dt\left( \tau  \right)}}{{d\tau }}d\tau }  = \frac{1}{\pi }\int_0^{ + \infty } {e^{ - \nu \tau } \frac{1}{{\sec \beta \cosh t\left( \tau  \right) + 1}}d\tau } .
\]
Following Howls, we express the function involving $t\left(\tau\right)$ as a contour integral using the residue theorem, to find
\[
\mathbf{A}_\nu  \left( {\nu \sec \beta } \right) = \frac{1}{\pi }\int_0^{ + \infty } {e^{ - \nu \tau } \frac{1}{{2\pi i}}\oint_\Gamma  {\frac{{f^{ - 1} \left( {u,\beta } \right)}}{{1 - \tau ^2 f^{ - 2} \left( {u,\beta } \right)}}du} d\tau } 
\]
where the contour $\Gamma$ encircles the path $\mathscr{P}\left(0\right)$ in the positive direction and does not enclose any of the saddle points $t_ \pm ^{\left( k \right)}$ (see Figure \ref{fig1}). Now, we employ the well-known expression for non-negative integer $N$
\begin{equation}\label{eq3}
\frac{1}{1 - z} = \sum\limits_{n = 0}^{N-1} {z^n}  + \frac{z^N}{1 - z},\; z \neq 1,
\end{equation}
to expand the function under the contour integral in powers of $\tau ^2 f^{ - 2} \left( {u,\beta } \right)$. The result is
\[
\mathbf{A}_\nu  \left( {\nu \sec \beta } \right) = \frac{1}{\pi}\sum\limits_{n = 0}^{N - 1}{ \int_0^{ + \infty } {\tau ^{2n} e^{ - \nu \tau } \frac{1}{{2\pi i}}\oint_\Gamma  {\frac{{du}}{{f^{2n + 1} \left( {u,\beta } \right)}}} d\tau }}  + R_N \left( {\nu ,\sec\beta } \right),
\]
where
\begin{equation}\label{eq4}
R_N \left( {\nu ,\sec\beta } \right) = \frac{1}{\pi }\int_0^{ + \infty } {\tau ^{2N} e^{ - \nu \tau } \frac{1}{{2\pi i}}\oint_\Gamma  {\frac{{f^{ - 2N - 1} \left( {u,\beta } \right)}}{{1 - \tau ^2 f^{ - 2} \left( {u,\beta } \right)}}du} d\tau } .
\end{equation}
The path $\Gamma$ in the sum can be shrunk into a small circle around $0$, and we arrive at
\begin{equation}\label{eq5}
\mathbf{A}_\nu  \left( {\nu \sec \beta } \right) = \frac{1}{\pi }\sum\limits_{n = 0}^{N - 1} {\frac{{\left( {2n} \right)!a_n \left( {\sec \beta } \right)}}{{\nu ^{2n + 1} }}}  + R_N \left( {\nu ,\sec\beta } \right),
\end{equation}
where
\[
a_n \left( {\sec \beta } \right) = \frac{1}{{2\pi i}}\oint_{\left( {0^ +  } \right)} {\frac{{du}}{{f^{2n + 1} \left( {u,\beta } \right)}}}  = \frac{1}{{\left( {2n} \right)!}}\left[ {\frac{{d^{2n} }}{{dt^{2n} }}\left( {\frac{t}{{\sec \beta \sinh t + t}}} \right)^{2n + 1} } \right]_{t = 0} .
\]

\begin{figure}[t]
\def\svgwidth{0.66\columnwidth}
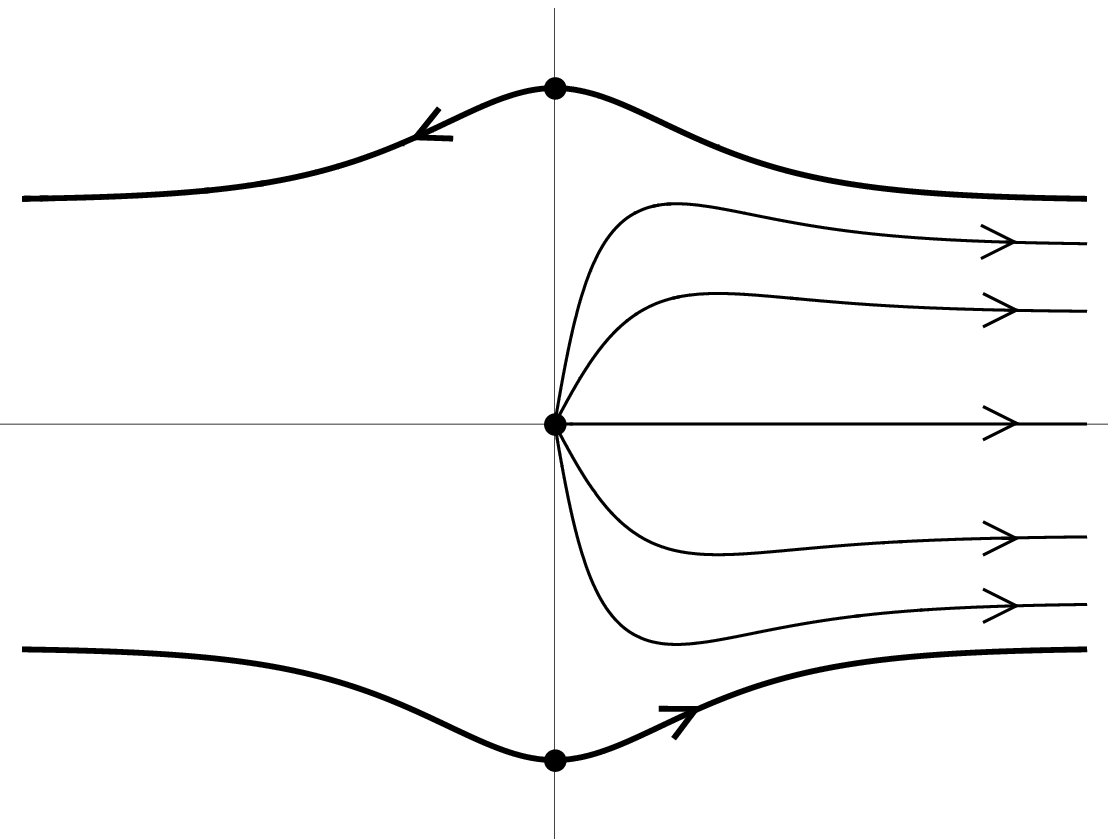
\caption{The path $\mathscr{P}\left(\theta\right)$ emanating from the origin when (i) $\theta=0$, (ii) $\theta=-\frac{\pi}{4}$, (iii) $\theta=-\frac{2\pi}{5}$, (iv) $\theta=\frac{\pi}{4}$, (v) $\theta=\frac{2\pi}{5}$. The paths $\mathscr{C}_{-} ^{\left( 0 \right)} \left( { -\frac{\pi }{2}} \right)$ and $\mathscr{C}_+ ^{\left( -1 \right)} \left( { \frac{\pi }{2}} \right)$ are the adjacent contours to $0$. The domain $\Delta$ comprises all points between these two paths in the right-half plane.}
\label{fig2}
\end{figure}

Performing the change of variable $\nu \tau = s$ in \eqref{eq4} yields
\begin{equation}\label{eq6}
R_N \left( {\nu ,\sec\beta } \right) = \frac{1}{{\pi \nu ^{2N + 1} }}\int_0^{ + \infty } {s^{2N} e^{ - s} \frac{1}{{2\pi i}}\oint_\Gamma  {\frac{{f^{ - 2N - 1} \left( {u,\beta } \right)}}{{1 - \left( {s/\nu } \right)^2 f^{ - 2} \left( {u,\beta } \right)}}du} ds} 
\end{equation}
This representation of $R_N \left( {\nu ,\sec\beta } \right)$ and the formula \eqref{eq5} can be continued analytically if we choose $\Gamma = \Gamma\left(\theta\right)$ to be an infinite contour that surrounds the path $\mathscr{P}\left(\theta\right)$ in the anti-clockwise direction and that does not encircle any of the saddle points $t_ \pm ^{\left( k \right)}$. This continuation argument works until the path $\mathscr{P}\left(\theta\right)$ runs into a saddle point. In the terminology of Howls, such saddle points are called adjacent to the endpoint $0$. As
\[
\left|\arg \left( {f\left(0 ,\beta\right) - f\left( {t_ \pm ^{\left( k \right)} ,\beta } \right)} \right)\right| = \frac{\pi}{2}
\]
for any saddle point $t_\pm ^{\left( k \right)}$, we infer that \eqref{eq6} is valid as long as $-\frac{\pi}{2} < \theta < \frac{\pi}{2}$ with a contour $\Gamma\left(\theta\right)$ specified above. When $\theta = -\frac{\pi}{2}$, the path $\mathscr{P}\left(\theta\right)$ connects to the saddle point $t_{-} ^{\left( 0 \right)} = - i\beta+\pi i$. Similarly, when $\theta = \frac{\pi}{2}$, the path $\mathscr{P}\left(\theta\right)$ connects to the saddle point $t_{+} ^{\left( -1 \right)} = i\beta-\pi i$. These are the adjacent saddles. The set
\[
\Delta  = \left\{ {u \in \mathscr{P}\left( \theta  \right) : - \frac{\pi}{2} < \theta  < \frac{\pi}{2}} \right\}
\]
forms a domain in the complex plane whose boundary contains portions of steepest descent paths through the adjacent saddles (see Figure \ref{fig2}). These paths are $\mathscr{C}_{-}^{\left( 0 \right)} \left( { -\frac{\pi }{2}} \right)$ and $\mathscr{C}_{+} ^{\left( -1 \right)} \left( {\frac{\pi }{2}} \right)$, and they are called the adjacent contours to the endpoint $0$. The function under the contour integral in \eqref{eq6} is an analytic function of $u$ in the domain $\Delta$, therefore we can deform $\Gamma$ over the adjacent contours. We thus find that for $-\frac{\pi}{2} < \theta < \frac{\pi}{2}$ and $N \geq 0$, \eqref{eq6} may be written
\begin{gather}\label{eq7}
\begin{split}
R_N \left( {\nu ,\sec\beta } \right) = \; & \frac{1}{{\pi \nu ^{2N + 1} }}\int_0^{ + \infty } {s^{2N} e^{ - s} \frac{1}{{2\pi i}}\int_{\mathscr{C}_{-}^{\left( 0 \right)} \left( { -\frac{\pi }{2}} \right)} {\frac{{f^{ - 2N - 1} \left( {u,\beta } \right)}}{{1 - \left( {s/\nu } \right)^2 f^{ - 2} \left( {u,\beta } \right)}}du} ds} \\ & + \frac{1}{{\pi \nu ^{2N + 1} }}\int_0^{ + \infty } {s^{2N} e^{ - s} \frac{1}{{2\pi i}}\int_{\mathscr{C}_{+} ^{\left( -1 \right)} \left( {\frac{\pi }{2}} \right)} {\frac{{f^{ - 2N - 1} \left( {u,\beta } \right)}}{{1 - \left( {s/\nu } \right)^2 f^{ - 2} \left( {u,\beta } \right)}}du} ds} .
\end{split}
\end{gather}

Now we make the changes of variable
\[
s = t\frac{{\left| {f\left( { - i\beta  + \pi i,\beta } \right) - f\left( {0,\beta } \right)} \right|}}{{f\left( { - i\beta  + \pi i,\beta } \right) - f\left( {0,\beta } \right)}}f\left( {u,\beta } \right) =  - i t f\left( {u,\beta } \right)
\]
in the first, and
\[
s = t\frac{{\left| {f\left( {i\beta  - \pi i,\beta } \right) - f\left( {0,\beta } \right)} \right|}}{{f\left( {i\beta  - \pi i,\beta } \right) - f\left( {0,\beta } \right)}}f\left( {u,\beta } \right) = i t f\left( {u,\beta } \right)
\]
in the second double integral. Clearly, by the definition of the adjacent contours, $t$ is positive. The quantities $f\left( { - i\beta  + \pi i,\beta } \right) - f\left( {0,\beta } \right) = i\left( {\tan \beta  - \beta +\pi} \right)$ and $f\left( {i\beta  - \pi i,\beta } \right) - f\left( {0,\beta } \right) = -i\left( {\tan \beta  - \beta +\pi } \right)$ were essentially called the ``singulants" by Dingle \cite[p. 147]{Dingle}. With these changes of variable, the representation \eqref{eq7} for $R_N \left( {\nu ,\sec\beta } \right)$ becomes
\begin{equation}\label{eq8}
R_N \left( {\nu ,\sec\beta } \right) = \frac{{\left( { - 1} \right)^N }}{{\pi \nu ^{2N + 1} }}\int_0^{ + \infty } {\frac{{t^{2N} }}{{1 + \left( {t/\nu } \right)^2 }}\left( {\frac{1}{{2\pi }}\int_{\mathscr{C}_{+} ^{\left( -1 \right)} \left( {\frac{\pi }{2}} \right)} {e^{ - itf\left( {u,\beta } \right)} du}  - \frac{1}{{2\pi }}\int_{\mathscr{C}_{-}^{\left( 0 \right)} \left( { -\frac{\pi }{2}} \right)} {e^{itf\left( {u,\beta } \right)} du} } \right)ds} ,
\end{equation}
for $-\frac{\pi}{2} < \theta < \frac{\pi}{2}$ and $N \geq 0$. Finally, we shift the contour $\mathscr{C}_{+} ^{\left( -1 \right)} \left( {\frac{\pi }{2}} \right)$ upward by $\pi$ and the contour $\mathscr{C}_{-}^{\left( 0 \right)} \left( { -\frac{\pi }{2}} \right)$ downward by $\pi$. Let us denote these new paths by $\widetilde{\mathscr{C}}_{+} ^{\left( -1 \right)} \left( {\frac{\pi }{2}} \right)$ and $\widetilde{\mathscr{C}}_{-}^{\left( 0 \right)} \left( { -\frac{\pi }{2}} \right)$, respectively. We therefore find that the contour integrals in \eqref{eq8} can be expressed in terms of the Hankel functions since
\[
\frac{1}{{2\pi }}\int_{\mathscr{C}_{+} ^{\left( -1 \right)} \left( {\frac{\pi }{2}} \right)} {e^{ - itf\left( {u,\beta } \right)} du}  = \frac{{e^{ - \pi t} }}{2}i\frac{1}{{\pi i}}\int_{\widetilde{\mathscr{C}}_{+} ^{\left( -1 \right)} \left( {\frac{\pi }{2}} \right)} {e^{it\left( {\sec \beta \sinh u - u} \right)} du}  = \frac{{e^{ - \pi t} }}{2}iH_{it}^{\left( 1 \right)} \left( {it\sec \beta } \right),
\]
and
\begin{align*}
 - \frac{1}{{2\pi }}\int_{\mathscr{C}_{-}^{\left( 0 \right)} \left( { -\frac{\pi }{2}} \right)} {e^{itf\left( {u,\beta } \right)} du}  =  - \frac{{e^{ - \pi t} }}{2}i\frac{1}{{\pi i}}\int_{\widetilde{\mathscr{C}}_{-}^{\left( 0 \right)} \left( { -\frac{\pi }{2}} \right)} {e^{ - it\left( {\sec \beta \sinh u - u} \right)} du}  & =  - \frac{{e^{ - \pi t} }}{2}iH_{ - it}^{\left( 2 \right)} \left( { - it\sec \beta } \right) \\ & = \frac{{e^{ - \pi t} }}{2}iH_{it}^{\left( 1 \right)} \left( {it\sec \beta } \right).
\end{align*}
Substituting these into \eqref{eq8} gives
\[
R_N \left( {\nu ,\sec\beta } \right) = \frac{{\left( { - 1} \right)^N }}{{\pi \nu ^{2N + 1} }}\int_0^{ + \infty } {\frac{{t^{2N} e^{ - \pi t} }}{{1 + \left( {t/\nu } \right)^2 }}iH_{it}^{\left( 1 \right)} \left( {it\sec \beta } \right)dt} ,
\]
for $-\frac{\pi}{2} < \theta < \frac{\pi}{2}$ and $N \geq 0$. Thus, we have proved \eqref{eq11} and \eqref{eq10} for $\lambda>1$.

Now, we extend our results to every $\lambda>0$. For fixed $\nu$, $\Re\left(\nu\right)>0$, we can extend $R_N \left( {\nu ,\lambda } \right)$ to an analytic function of $\lambda>0$ using \eqref{eq11}. From the known behaviours
\[
iH_{it}^{\left( 1 \right)} \left( {\lambda it} \right) \sim  - \frac{2}{\pi }\log t \; \text{ as } \; t \to 0+
\]
and
\[
iH_{it}^{\left( 1 \right)} \left( {\lambda it} \right) = o\left( 1 \right) \; \text{ as } \; t \to +\infty,
\]
it is seen that the right-hand side of \eqref{eq10} is a well-defined analytic function of $0<\lambda \leq 1$, for every fixed $\nu$ with $\Re\left(\nu\right)>0$. Whence, by analytic continuation the equality \eqref{eq10} holds for every $\lambda>0$ and $\nu$ with $\Re\left(\nu\right)>0$.

The first formula in \eqref{eq9} has been proved for $\lambda>1$, however, by analytic continuation, it holds for every $\lambda>0$. To prove the second representation in \eqref{eq9}, we apply \eqref{eq10} for the right-hand side of
\[
a_n \left( \lambda  \right) = \pi \frac{{\nu ^{2n + 1} }}{{\left( {2n} \right)!}}\left( {R_n \left( {\nu ,\lambda } \right) - R_{n + 1} \left( {\nu ,\lambda } \right)} \right).
\]

\section{Error bounds}\label{section3}

In this section, we derive explicit and realistic error bounds for the large order asymptotic series of the Anger--Weber function. The proofs are based on the resurgence formula given in Theorem \ref{thm1}.

We comment on the relation between Meijer's work \cite{Meijer} on the asymptotic expansion of $\mathbf{A}_{\nu} \left( \lambda \nu \right)$, $\lambda > 1$ and ours. Some of the estimates in \cite{Meijer} coincide with ours and are valid in wider sectors of the complex $\nu$-plane. However, it should be noted that those bounds become less effective outside the sectors of validity of the representation \eqref{eq10} due to the Stokes phenomenon. For those sectors we recommend the use of the continuation formulas given in Section \ref{section1}.

To estimate the remainder terms, we shall use the elementary result that
\begin{equation}\label{eq12}
\frac{1}{{\left| {1 - re^{i\varphi } } \right|}} \le \begin{cases} \left|\csc \varphi \right| & \; \text{ if } \; 0 < \left|\varphi \text{ mod } 2\pi\right| <\frac{\pi}{2} \\ 1 & \; \text{ if } \; \frac{\pi}{2} \leq \left|\varphi \text{ mod } 2\pi\right| \leq \pi \end{cases}
\end{equation}
holds for any $r>0$. We will also need the fact that
\begin{equation}\label{eq13}
iH_{it}^{\left(1\right)} \left( {\lambda it} \right) \ge 0
\end{equation}
for any $t>0$ and $\lambda \geq 1$ (see \cite{Nemes}).

\subsection{Case (i): $\lambda \geq 1$} We observe that from \eqref{eq9} and \eqref{eq13} it follows that
\[
\left| {a_n \left( \lambda  \right)} \right| = \frac{1}{{\left( {2n} \right)!}}\int_0^{ + \infty } {t^{2n} e^{ - \pi t} iH_{it}^{\left( 1 \right)} \left( {\lambda it} \right)dt} .
\]
Using this formula, together with the representation \eqref{eq10} and the estimate \eqref{eq12}, we obtain the error bound
\begin{equation}\label{eq14}
\left| {R_N \left( {\nu ,\lambda } \right)} \right| \le \frac{1}{\pi }\frac{{\left( {2N} \right)!\left| {a_N \left( \lambda  \right)} \right|}}{{\left| \nu  \right|^{2N + 1} }} \begin{cases} \left|\csc\left(2\theta\right)\right| & \; \text{ if } \; \frac{\pi}{4} < \left|\theta\right| <\frac{\pi}{2} \\ 1 & \; \text{ if } \; \left|\theta\right| \leq \frac{\pi}{4}. \end{cases}
\end{equation}
When $\nu$ is real and positive, we can obtain more precise estimates. Indeed, as $0 < \frac{1}{{1 + \left( {t/\nu } \right)^2 }} < 1 $ for $t,\nu>0$, from \eqref{eq10} and \eqref{eq9} we find
\[
R_N \left( {\nu ,\lambda } \right) = \frac{1}{\pi }\frac{{\left( {2N} \right)!a_N \left( \lambda  \right)}}{{\nu ^{2N + 1} }}\Theta ,
\]
where $0 < \Theta < 1$ is an appropriate number depending on $\nu,\lambda$ and $N$. In particular, when $N=0$, we have
\[
0 < \mathbf{A}_\nu  \left( {\lambda \nu } \right) < \frac{1}{{\pi \nu \left( {1 + \lambda } \right)}} \; \text{ for } \; \nu >0.
\]
Therefore, the leading order asymptotic approximation for $\mathbf{A}_{\nu} \left( \lambda \nu \right)$ is always in error by excess, for $\lambda \geq 1$ and for all positive values of $\nu$.

The error bound \eqref{eq14} becomes singular as $\theta \to \pm \frac{\pi}{2}$, and therefore unrealistic near the Stokes lines. A better bound for $R_N \left( {\nu ,\lambda } \right)$ near these lines can be derived as follows. Let $0 < \varphi  < \frac{\pi }{2}$ be an acute angle that may depend on $N$. Suppose that $\frac{\pi}{4} +\varphi < \theta \le \frac{\pi}{2}$. An analytic continuation of the representation \eqref{eq11} to this sector can be found by rotating the path of integration in \eqref{eq10} by $\varphi$:
\[
R_N \left( {\nu ,\lambda } \right) = \frac{{\left( { - 1} \right)^N }}{{\pi \nu ^{2N + 1} }}\int_0^{ + \infty e^{i\varphi } } {\frac{{t^{2N} e^{ - \pi t} }}{{1 + \left( {t/\nu } \right)^2 }}iH_{it}^{\left( 1 \right)} \left( {\lambda it} \right)dt} .
\]
Substituting $t = \frac{se^{i\varphi }}{\cos \varphi}$ and applying the estimation \eqref{eq12}, we obtain
\[
\left| {R_N \left( {\nu ,\lambda } \right)} \right| \le \frac{{\csc \left( {2\left( {\theta  - \varphi } \right)} \right)}}{{\pi \cos ^{2N + 1} \varphi \left| \nu  \right|^{2N + 1} }}\int_0^{ + \infty } {s^{2N} e^{ - \pi s} \left| {H_{\frac{{ise^{i\varphi } }}{{\cos \varphi }}}^{\left( 1 \right)} \left( {\lambda \frac{{ise^{i\varphi } }}{{\cos \varphi }}} \right)} \right|ds} .
\]
In \cite{Nemes}, it was shown that
\begin{equation}\label{eq60}
\left| {H_{\frac{{ise^{i\varphi } }}{{\cos \varphi }}}^{\left( 1 \right)} \left( {\lambda \frac{{ise^{i\varphi } }}{{\cos \varphi }}} \right)} \right| \le \frac{1}{{\sqrt {\cos \varphi } }}\left| {H_{is}^{\left( 1 \right)} \left( {\lambda is} \right)} \right| = \frac{1}{{\sqrt {\cos \varphi } }}iH_{is}^{\left( 1 \right)} \left( {\lambda is} \right)
\end{equation}
for any $s>0$, $\lambda\geq 1$ and $0 < \varphi  < \frac{\pi }{2}$. It follows that
\begin{equation}\label{eq15}
\left| {R_N \left( {\nu ,\lambda } \right)} \right| \le \frac{{\csc \left( {2\left( {\theta  - \varphi } \right)} \right)}}{{\pi \cos ^{2N + \frac{3}{2}} \varphi \left| \nu  \right|^{2N + 1} }}\int_0^{ + \infty } {s^{2N} e^{ - \pi s} iH_{is}^{\left( 1 \right)} \left( {\lambda is} \right)ds}  = \frac{{\csc \left( {2\left( {\theta  - \varphi } \right)} \right)}}{{\cos ^{2N + \frac{3}{2}} \varphi }}\frac{1}{\pi }\frac{{\left( {2N} \right)!\left| {a_N \left( \lambda  \right)} \right|}}{{\left| \nu  \right|^{2N + 1} }}.
\end{equation}
The angle $\varphi  = \arctan \left( {\left( {\frac{{4N + 5}}{2}} \right)^{ - \frac{1}{2}} } \right)$ minimizes the function $\csc \left( {2\left( {\frac{\pi }{2} - \varphi } \right)} \right)\cos ^{ - 2N - \frac{3}{2}} \varphi$, and
\begin{align*}
\frac{{\csc \left( {2\left( {\theta  - \arctan \left( {\left( {\frac{{4N + 5}}{2}} \right)^{ - \frac{1}{2}} } \right)} \right)} \right)}}{{\cos ^{2N + \frac{3}{2}} \left( {\arctan \left( {\left( {\frac{{4N + 5}}{2}} \right)^{ - \frac{1}{2}} } \right)} \right)}} & \le \frac{{\csc \left( {2\left( {\frac{\pi }{2} - \arctan \left( {\left( {\frac{{4N + 5}}{2}} \right)^{ - \frac{1}{2}} } \right)} \right)} \right)}}{{\cos ^{2N + \frac{3}{2}} \left( {\arctan \left( {\left( {\frac{{4N + 5}}{2}} \right)^{ - \frac{1}{2}} } \right)} \right)}} \\ &= \frac{1}{{\sqrt 2 }}\left( 1+\frac{2}{4N + 5} \right)^{N + \frac{7}{4}} \sqrt {N + \frac{5}{4}}  \le \sqrt {\frac{e}{2}\left( {N + \frac{3}{2}} \right)}
\end{align*}
for all $\frac{\pi }{4} + \varphi  = \frac{\pi }{4} + \arctan \left( {\left( {\frac{{4N + 5}}{2}} \right)^{ - \frac{1}{2}} } \right) < \theta  \le \frac{\pi }{2}$ with $N \geq 0$. Applying this in \eqref{eq15} yields the upper bound
\begin{equation}\label{eq16}
\left| {R_N \left( {\nu ,\lambda } \right)} \right| \le \sqrt {\frac{e}{2}\left( {N + \frac{3}{2}} \right)} \frac{1}{\pi }\frac{{\left( {2N} \right)!\left| {a_N \left( \lambda  \right)} \right|}}{{\left| \nu  \right|^{2N + 1} }},
\end{equation}
which is valid for $\frac{\pi }{4} + \varphi  = \frac{\pi }{4} + \arctan \left( {\left( {\frac{{4N + 5}}{2}} \right)^{ - \frac{1}{2}} } \right) < \theta  \le \frac{\pi }{2}$ with $N \geq 0$. Since $\left| {R_N \left( {\bar \nu ,\lambda } \right)} \right| = \left| {\overline {R_N \left( {\nu ,\lambda } \right)} } \right| = \left| {R_N \left( {\nu ,\lambda } \right)} \right|$, this bound also holds when $-\frac{\pi}{2} \leq \theta < -\frac{\pi }{4} - \arctan \left( {\left( {\frac{{4N + 5}}{2}} \right)^{ - \frac{1}{2}} } \right)$. In the ranges $\frac{\pi }{4} < \left| \theta  \right| \leq \frac{\pi }{4} + \arctan \left( {\frac{{\sqrt 2 }}{3}} \right)$ it holds that $\left| {\csc \left( {2\theta } \right)} \right| \le \sqrt {\frac{e}{2}\left( {1 + \frac{3}{2}} \right)}$, whence the estimate \eqref{eq16} is valid in the wider sectors $\frac{\pi }{4} < \left| \theta  \right| \le \frac{\pi }{2}$ as long as $N\geq 1$.

\subsection{Case (ii): $0 < \lambda < 1$} In this case, we cannot prove error bounds involving the first omitted term, since $iH_{it}^{\left(1\right)} \left( {\lambda it} \right)$ has an oscillatory behaviour when $0 < \lambda < 1$ and $t>0$. Nevertheless, we define $\widetilde a_n \left( \lambda  \right)$ via the integral
\[
\widetilde a_n \left( \lambda  \right) = \frac{1}{{\left( {2n} \right)!}}\int_0^{ + \infty } {t^{2n} e^{ - \pi t} \left| {H_{it}^{\left( 1 \right)} \left( {\lambda it} \right)} \right|dt};
\]
and by the representation \eqref{eq10} and the inequality \eqref{eq12}, we deduce the error bound
\[
\left| {R_N \left( {\nu ,\lambda } \right)} \right| \le \frac{1}{\pi }\frac{{\left( {2N} \right)! \widetilde a_N \left( \lambda  \right)}}{{\left| \nu  \right|^{2N + 1} }} \begin{cases} \left|\csc\left(2\theta\right)\right| & \; \text{ if } \; \frac{\pi}{4} < \left|\theta\right| <\frac{\pi}{2} \\ 1 & \; \text{ if } \; \left|\theta\right| \leq \frac{\pi}{4}. \end{cases}
\]
Simple estimates for the quantities $\widetilde a_n \left(\lambda\right)$ may perhaps be derived from the connection formula with the modified Bessel function of the third kind of purely imaginary order
\[
H_{it}^{\left( 1 \right)} \left(\lambda it \right) = \frac{2}{\pi i}e^{\frac{\pi}{2}t} K_{it} \left( \lambda t \right),
\]
and the known bounds for this latter function (see, e.g., Booker et al. \cite{Booker}).

Since for $0<\lambda <1$ we do not have an inequality like \eqref{eq60}, it seems hard to obtain any usable error bound which is appropriate when $\arg \nu$ is close to $\pm \frac{\pi}{2}$.

\section{Asymptotics for the late coefficients}\label{section4}

In this section, we investigate the asymptotic nature of the coefficients $a_n\left(\lambda\right)$ as $n \to +\infty$ with $\lambda$ being fixed. For our purposes, the most appropriate representation of these coefficients is the second integral formula in \eqref{eq9}. Although the representation is valid for all $\lambda>0$, we shall find that the asymptotic form of $a_n\left(\lambda\right)$ is significantly different according to whether $\lambda>1$, $\lambda=1$ or $0<\lambda<1$.

\subsection{Case (i): $\lambda>1$} For this case, we take $\lambda = \sec \beta$ with a suitable $0 < \beta <\frac{\pi}{2}$. From \eqref{eq30}, it follows that for any $t>0$ and $0 < \beta <\frac{\pi}{2}$, it holds that
\begin{equation}\label{eq17}
iH_{it}^{\left( 1 \right)} \left( {it\sec \beta } \right) = \frac{{e^{ - t\left( {\tan \beta  - \beta } \right)} }}{{\left( {\frac{1}{2}t\pi \tan \beta } \right)^{\frac{1}{2}} }}\left( {\sum\limits_{m = 0}^{M - 1} {\frac{{i^m U_m \left( {i\cot \beta } \right)}}{{t^m }}}  + R_M^{\left( H \right)} \left( {it,\beta } \right)} \right).
\end{equation}
In \cite{Nemes}, it was proved that the remainder $R_M^{\left( H \right)} \left( {it,\beta } \right)$ satisfies
\begin{equation}\label{eq18}
\left| {R_M^{\left( H \right)} \left( {it,\beta } \right)} \right| \le \frac{{\left| {U_M \left( {i\cot \beta } \right)} \right|}}{{t^M }} .
\end{equation}
Substituting the formula \eqref{eq17} into \eqref{eq9} leads us to the expansion
\begin{gather}\label{eq19}
\begin{split}
\left( {2n} \right)!a_n \left( {\sec \beta } \right) = \; & \left( {\frac{{2\cot \beta }}{{\pi \left( {\tan \beta  - \beta  + \pi } \right)}}} \right)^{\frac{1}{2}} \frac{{\left( { - 1} \right)^n \Gamma \left( {2n + \frac{1}{2}} \right)}}{{\left( {\tan \beta  - \beta  + \pi } \right)^{2n} }} \\ & \times \left( {\sum\limits_{m = 0}^{M - 1} {\left( {i\left( {\tan \beta  - \beta  + \pi } \right)} \right)^m U_m \left( {i\cot \beta } \right)\frac{{\Gamma \left( {2n - m + \frac{1}{2}} \right)}}{{\Gamma \left( {2n + \frac{1}{2}} \right)}}}  + A_M \left( {n,\beta } \right)} \right),
\end{split}
\end{gather}
for any fixed $0 \le M \le 2n$, provided that $n\geq 1$. The remainder term $A_M \left( {n,\beta } \right)$ is given by the integral formula
\[
A_M \left( {n,\beta } \right) = \frac{{\left( {\tan \beta  - \beta  + \pi } \right)^{2n + \frac{1}{2}} }}{{\Gamma \left( {2n + \frac{1}{2}} \right)}}\int_0^{ + \infty } {t^{2n - \frac{1}{2}} e^{ - \left( {\tan \beta  - \beta  + \pi } \right)t} R_M^{\left( H \right)} \left( {it,\beta } \right)dt} .
\]
To bound this error term, we apply the estimate \eqref{eq18} to find
\begin{equation}\label{eq20}
\left| {A_M \left( {n,\beta } \right)} \right| \le \left( {\tan \beta  - \beta  + \pi } \right)^M \left| {U_M \left( {i\cot \beta } \right)} \right|\frac{{\Gamma \left( {2n - M + \frac{1}{2}} \right)}}{{\Gamma \left( {2n + \frac{1}{2}} \right)}}.
\end{equation}
Expansions of type \eqref{eq19} are called inverse factorial series in the literature. Numerically, their character is similar to the character of asymptotic power series, because the consecutive Gamma functions decrease asymptotically by a factor $2n$.

From the asymptotic behaviour of the coefficients $U_m \left( {i\cot \beta } \right)$ (see \cite{Nemes}), we infer that for large $n$, the least value of the bound \eqref{eq20} occurs when
\[
M \approx \frac{\tan \beta  - \beta }{3\left( {\tan \beta  - \beta } \right) + \pi}\left(4n+1\right) .
\]
Whence, the smaller $\beta$ is the larger $n$ has to be to get a reasonable approximation from \eqref{eq19}.

Numerical examples illustrating the efficacy of the expansion \eqref{eq19}, truncated optimally, are given in Table \ref{table1}.

\begin{table*}[!ht]
\begin{center}
\begin{tabular}
[c]{ l r @{\,}c@{\,} l}\hline
 & \\ [-1ex]
 values of $\beta$ and $M$ & $\beta=\frac{\pi}{6}$, $M=4$ & & \\ [1ex]
 exact numerical value of $a_{50}\left(\sec\beta\right)$ & $0.2004926124399177097019512509947129$ & $\times$ & $10^{-51}$ \\ [1ex]
 approximation \eqref{eq19} to $a_{50}\left(\sec\beta\right)$ & $0.1997204566354320191164985775448290$ & $\times$ & $10^{-51}$  \\ [1ex]
 error & $0.7721558044856905854526734498839$ & $\times$ & $10^{-54}$\\ [1ex] 
 error bound using \eqref{eq20} & $0.16182537012652011778281419657176$ & $\times$ & $10^{-53}$\\ [1ex] \hline
 & \\ [-1ex]
 values of $\beta$ and $M$ &  $\beta=\frac{\pi}{3}$, $M=27$ & & \\ [1ex]
 exact numerical value of $a_{50}\left(\sec\beta\right)$ & $0.1619316740481494064448396260188866$ & $\times$ & $10^{-59}$ \\ [1ex]
  approximation \eqref{eq19} to $a_{50}\left(\sec\beta\right)$ & $0.1619316740481497277978573226174596$ & $\times$ & $10^{-59}$ \\ [1ex]
 error & $-0.3213530176965985730$ & $\times$ & $10^{-74}$\\ [1ex]
 error bound using \eqref{eq20} & $0.6473043619300051742$ & $\times$ & $10^{-74}$\\ [1ex] \hline
 & \\ [-1ex]
 values of $\beta$ and $M$ &  $\beta=\frac{5\pi}{12}$, $M=47$ & & \\ [1ex]
 exact numerical value of $a_{50}\left(\sec\beta\right)$ & $0.4989354184460076118014557886550703$ & $\times$ & $10^{-76}$ \\ [1ex]
 approximation \eqref{eq19} to $a_{50}\left(\sec\beta\right)$ & $0.4989354184460076118014557886641359$ & $\times$ & $10^{-76}$ \\ [1ex]
 error & $-0.90656$ & $\times$ & $10^{-105}$\\ [1ex] 
 error bound using \eqref{eq20} & $0.181989$ & $\times$ & $10^{-104}$\\ [-1ex]
 & \\\hline
\end{tabular}
\end{center}
\caption{Approximations for $a_{50}\left(\sec\beta\right)$ with various $\beta$, using \eqref{eq19}.}
\label{table1}
\end{table*}

\subsection{Case (ii): $\lambda=1$} Using \eqref{eq31}, we can write
\begin{equation}\label{eq24}
iH_{it}^{\left( 1 \right)} \left( {it} \right) = \frac{2}{{3\pi }}\sum\limits_{m = 0}^{M - 1} {\left( { - 1} \right)^m d_{2m} \sin \left( {\frac{{\left( {2m + 1} \right)\pi }}{3}} \right)\frac{{\Gamma \left( {\frac{{2m + 1}}{3}} \right)}}{{t^{\frac{{2m + 1}}{3}} }}}  + iR_M^{\left( H \right)} \left( {it} \right)
\end{equation}
for any $t>0$. It was shown in \cite{Nemes} that
\begin{equation}\label{eq25}
\left| {R_M^{\left( H \right)} \left( {it} \right)} \right| \le \frac{2}{3\pi}\left| {d_{2M} } \right|\frac{\sqrt 3 }{2}\frac{{\Gamma \left( {\frac{2M + 1}{3}} \right)}}{{t^{\frac{2M + 1}{3}} }} + \frac{2}{3\pi}\left| {d_{2M + 2} } \right|\frac{\sqrt 3}{2}\frac{{\Gamma \left( {\frac{2M + 3}{3}} \right)}}{{t^{\frac{2M + 3}{3}} }},
\end{equation}
\begin{equation}\label{eq26}
\left| {R_M^{\left( H \right)} \left( {it} \right)} \right| \le \frac{2}{3\pi}\left| {d_{2M + 2} } \right|\frac{\sqrt 3}{2}\frac{{\Gamma \left( {\frac{2M + 3}{3}} \right)}}{{t^{\frac{2M + 3}{3}} }},
\end{equation}
\begin{equation}\label{eq27}
\left| {R_M^{\left( H \right)} \left( {it} \right)} \right| \le \frac{2}{3\pi}\left| {d_{2M} } \right|\frac{\sqrt 3}{2}\frac{{\Gamma \left( {\frac{2M + 1}{3}} \right)}}{{t^{\frac{2M + 1}{3}} }}
\end{equation}
according to whether $M \equiv 0 \mod 3$, $M \equiv 1 \mod 3$ or $M \equiv 2 \mod 3$, respectively. Substituting the expression \eqref{eq24} into \eqref{eq9} yields the expansion
\begin{gather}\label{eq21}
\begin{split}
\left( {2n} \right)!a_n \left( 1 \right) = \; & \left( { - 1} \right)^n \frac{2\Gamma \left( {2n + \frac{2}{3}} \right)}{3\pi ^{2n + \frac{5}{3}} }\\ & \times \left( {\sum\limits_{m = 0}^{M - 1} {\left( { - 1} \right)^m \pi ^{\frac{{2m}}{3}} d_{2m} \sin \left( {\frac{{\left( {2m + 1} \right)\pi }}{3}} \right)\Gamma \left( {\frac{{2m + 1}}{3}} \right)\frac{{\Gamma \left( {\frac{{6n - 2m + 2}}{3}} \right)}}{{\Gamma \left( {2n + \frac{2}{3}} \right)}}}  + A_M \left( n \right)} \right),
\end{split}
\end{gather}
for any fixed $0 \leq M \leq 3n - 1$, provided that $n \geq 1$. The remainder term $A_M\left( n \right)$ is given by the formula
\[
A_M \left( n \right) = \frac{{3\pi ^{2n + \frac{5}{3}} }}{{2\Gamma \left( {2n + \frac{2}{3}} \right)}}\int_0^{ + \infty } {t^{2n} e^{ - \pi t} iR_M^{\left( H \right)} \left( {it} \right)dt} .
\]
Bounds for this error term follow from the estimates \eqref{eq25}--\eqref{eq27} since
\[
\left| {A_M \left( n \right)} \right| \le \pi ^{\frac{{2M}}{3}} \left| {d_{2M} } \right|\frac{{\sqrt 3 }}{2}\Gamma \left( {\frac{{2M + 1}}{3}} \right)\frac{{\Gamma \left( {\frac{{6n - 2M + 2}}{3}} \right)}}{{\Gamma \left( {2n + \frac{2}{3}} \right)}} + \pi ^{\frac{{2M + 2}}{3}} \left| {d_{2M + 2} } \right|\frac{{\sqrt 3 }}{2}\Gamma \left( {\frac{{2M + 3}}{3}} \right)\frac{{\Gamma \left( {\frac{{6n - 2M}}{3}} \right)}}{{\Gamma \left( {2n + \frac{2}{3}} \right)}},
\]
\begin{equation}\label{eq22}
\left| {A_M \left( n \right)} \right| \le \pi ^{\frac{{2M + 2}}{3}} \left| {d_{2M + 2} } \right|\frac{{\sqrt 3 }}{2}\Gamma \left( {\frac{{2M + 3}}{3}} \right)\frac{{\Gamma \left( {\frac{{6n - 2M}}{3}} \right)}}{{\Gamma \left( {2n + \frac{2}{3}} \right)}},
\end{equation}
\begin{equation}\label{eq23}
\left| {A_M \left( n \right)} \right| \le \pi ^{\frac{{2M}}{3}} \left| {d_{2M} } \right|\frac{{\sqrt 3 }}{2}\Gamma \left( {\frac{{2M + 1}}{3}} \right)\frac{{\Gamma \left( {\frac{{6n - 2M + 2}}{3}} \right)}}{{\Gamma \left( {2n + \frac{2}{3}} \right)}}
\end{equation}
according to whether $M\equiv 0 \mod 3$, $M\equiv 1 \mod 3$ or $M\equiv 2 \mod 3$, respectively.

From the asymptotic behaviour of the coefficients $d_{2m}$ (see \cite{Nemes}), for large $n$, the least values of these bounds occur when $M \approx 2n$. With this choice of $M$, the error bounds are $\mathcal{O}\left( {n^{ - \frac{1}{2}} 9^{-n} } \right)$. This is the best accuracy we can achieve using the expansion \eqref{eq21}. Numerical examples for various $n$ are provided in Table \ref{table2}.

\begin{table*}[!ht]
\begin{center}
\begin{tabular}
[c]{ l r @{\,}c@{\,} l}\hline
 & \\ [-1ex]
 values of $n$ and $M$ & $n=5$, $M=10$ & & \\ [1ex]
 exact numerical value of $a_n\left(1\right)$ & $-0.2039315629047481261022927689594356$ & $\times$ & $10^{-5}$ \\ [1ex]
 approximation \eqref{eq21} to $a_n\left(1\right)$ & $-0.2039317236866484733447636037370858$ & $\times$ & $10^{-5}$  \\ [1ex]
 error & $0.1607819003472424708347776502$ & $\times$ & $10^{-11}$\\ [1ex]
 error bound using \eqref{eq22} & $0.5218454726884724646870658288$ & $\times$ & $10^{-11}$\\ [1ex] \hline
 & \\ [-1ex]
 values of $n$ and $M$ & $n=10$, $M=20$ & & \\ [1ex]
 exact numerical value of $a_n\left(1\right)$ & $0.1740499192613222665759959822566006$ & $\times$ & $10^{-10}$ \\ [1ex]
 approximation \eqref{eq21} to $a_n\left(1\right)$ & $0.1740499192631695872689300620308834$ & $\times$ & $10^{-10}$  \\ [1ex]
 error & $-0.18473206929340797742828$ & $\times$ & $10^{-21}$\\ [1ex]
 error bound using \eqref{eq23} & $0.52455141471539645254342$ & $\times$ & $10^{-21}$\\ [1ex] \hline
 & \\ [-1ex]
 values of $n$ and $M$ & $n=25$, $M=50$ & &\\ [1ex]
 exact numerical value of $a_n\left(1\right)$ & $-0.1567780710784896492198553870128892$ & $\times$ & $10^{-25}$ \\ [1ex]
 approximation \eqref{eq21} to $a_n\left(1\right)$ & $-0.1567780710784896492198553919627602$ & $\times$ & $10^{-25}$ \\ [1ex]
 error & $0.49498710$ & $\times$ & $10^{-51}$\\ [1ex]
 error bound using \eqref{eq23} & $0.145150293$ & $\times$ & $10^{-50}$\\ [-1ex]
 & \\\hline
\end{tabular}
\end{center}
\caption{Approximations for $a_n\left(1\right)$ with various $n$, using \eqref{eq21}.}
\label{table2}
\end{table*}

\subsection{Case (iii): $0<\lambda<1$} For this case, we take $\lambda = \mathop{\text{sech}} \alpha$ with a suitable $\alpha>0$. It is known that
\begin{equation}\label{eq28}
iH_{it}^{\left( 1 \right)} \left( {it \mathop{\text{sech}} \alpha } \right) = 2 \Re \left( {\frac{{e^{it\left( {\alpha  - \tanh \alpha } \right) - \frac{\pi }{4}i} }}{{\left( {\frac{1}{2}\pi t\tanh \alpha } \right)^{\frac{1}{2}} }}\left( {\sum\limits_{m = 0}^{M - 1} {\frac{{i^m U_m \left( {\coth \alpha } \right)}}{{t^m }}}  + R_M^{\left( H \right)} \left( {it,\alpha } \right)} \right)} \right),
\end{equation}
where $R_M^{\left( H \right)} \left( {it,\alpha } \right) = \mathcal{O}_{M,\alpha } \left( {t^{ - M} } \right)$ as $t \to +\infty$. Here $U_m \left( {\coth \alpha } \right) = \left[U_m \left( x \right)\right]_{x = \coth \alpha }$ with $U_m\left(x\right)$ being a polynomial in $x$ of degree $3m$. As far as we know, there is no simple closed expression nor a realistic estimate for the remainder term $R_M^{\left( H \right)} \left( {it,\alpha } \right)$. Nevertheless, we assume that
\[
\int_0^{ + \infty } {t^{2n - \frac{1}{2}} e^{ - \pi t} \left| {R_M^{\left( H \right)} \left( {it,\alpha } \right)} \right|dt} < +\infty
\]
and substitute the expansion \eqref{eq28} into \eqref{eq9} to obtain
\begin{gather}\label{eq29}
\begin{split}
\left( {2n} \right)!a_n \left( {\mathop{\text{sech}} \alpha } \right) = \;& \Re \left( {\left( {\frac{{2\coth \alpha }}{{\pi \left( {\alpha  - \tanh \alpha  + \pi i} \right)}}} \right)^{\frac{1}{2}} \frac{2 \Gamma \left( {2n + \frac{1}{2}} \right)}{{\left( {\alpha  - \tanh \alpha  + \pi i} \right)^{2n} }} }\right.  \\ & \times \left.{\left( {\sum\limits_{m = 0}^{M - 1} {\left( {\alpha  - \tanh \alpha  + \pi i} \right)^m U_m \left( {\coth \alpha } \right)\frac{\Gamma \left( {2n - m + \frac{1}{2}} \right)}{\Gamma \left( {2n + \frac{1}{2}} \right)}}  + A_M \left( {n,\alpha } \right)} \right)} \right),
\end{split}
\end{gather}
for any fixed $0 \le M \le 2n$, provided that $n\geq 1$. The remainder term $A_M \left( {n,\alpha } \right)$ is given by the integral formula
\[
A_M \left( {n,\alpha } \right) = \left( { - 1} \right)^n e^{ - \frac{\pi }{4}i} \frac{{\left( {\alpha  - \tanh \alpha  + \pi i} \right)^{2n + \frac{1}{2}} }}{{\Gamma \left( {2n + \frac{1}{2}} \right)}}\int_0^{ + \infty } {t^{2n - \frac{1}{2}} e^{i\left( {\alpha  - \tanh \alpha  + \pi i} \right)t} R_M^{\left( H \right)} \left( {it,\alpha } \right)dt} .
\]
To achieve the best accuracy using the expansion \eqref{eq29}, we need to determine the index of the least term of the expansion. This can be done if we know the large $m$ behaviour of the coefficients $U_m \left(\coth \alpha \right)$. Such an asymptotic formula was derived by Dingle \cite[p. 168]{Dingle}, using formal, non-rigorous methods. At leading order, his formula can be written as
\[
\left| {U_m \left( {\coth \alpha } \right)} \right| \sim \frac{{\Gamma \left( m \right)}}{{2\pi \left( {2\left( {\alpha  - \tanh \alpha } \right)} \right)^m }}.
\]
Numerical calculations indicate that this approximation is correct, and assuming so, the optimal truncation occurs at
\[
M \approx \frac{{\alpha  - \tanh \alpha }}{{2\left( {\alpha  - \tanh \alpha } \right) + \sqrt {\left( {\alpha  - \tanh \alpha } \right)^2  + \pi ^2 } }}\left( {4n + 1} \right).
\]
Therefore, the smaller $\alpha$ is the larger $n$ has to be to get a reasonable approximation from \eqref{eq29}.

Numerical examples illustrating the efficacy of the formula \eqref{eq29}, truncated optimally, are given in Table \ref{table3}.

\begin{table*}[!ht]
\begin{center}
\begin{tabular}
[c]{ l r @{\,}c@{\,} l}\hline
 & \\ [-1ex]
 values of $\alpha$ and $M$ & $\alpha=\frac{1}{2}$, $M=3$ & & \\ [1ex]
 exact numerical value of $a_{50}\left(\mathop{\text{sech}} \alpha\right)$ & $0.2315627683882018769175712540082165$ & $\times$ & $10^{-50}$ \\ [1ex]
 approximation \eqref{eq29} to $a_{50}\left(\mathop{\text{sech}} \alpha\right)$ & $0.2303064844873166640986637287015961$ & $\times$ & $10^{-50}$  \\ [1ex]
 error & $0.12562839008852128189075253066203$ & $\times$ & $10^{-52}$\\ [1ex] \hline
 & \\ [-1ex]
 values of $\alpha$ and $M$ &  $\alpha=1$, $M=14$ & & \\ [1ex]
 exact numerical value of $a_{50}\left(\mathop{\text{sech}} \alpha\right)$ & $0.1279482878426982457824386451759845$ & $\times$ & $10^{-50}$ \\ [1ex]
  approximation \eqref{eq29} to $a_{50}\left(\mathop{\text{sech}} \alpha\right)$ & $0.1279482903067682761677730364825915$ & $\times$ & $10^{-50}$ \\ [1ex]
 error & $-0.24640700303853343913066070$ & $\times$ & $10^{-58}$\\ [1ex] \hline
 & \\ [-1ex]
 values of $\alpha$ and $M$ &  $\alpha=5$, $M=62$ & & \\ [1ex]
 exact numerical value of $a_{50}\left(\mathop{\text{sech}} \alpha\right)$ & $-0.9536145099812834565097014294181624$ & $\times$ & $10^{-72}$ \\ [1ex]
 approximation \eqref{eq29} to $a_{50}\left(\mathop{\text{sech}} \alpha\right)$ & $-0.9536145099812834565097014294180344$ & $\times$ & $10^{-72}$ \\ [1ex]
 error & $0.1280$ & $\times$ & $10^{-102}$\\ [-1ex] 
  & \\\hline
\end{tabular}
\end{center}
\caption{Approximations for $a_{50}\left(\mathop{\text{sech}} \alpha\right)$ with various $\alpha$, using \eqref{eq29}.}
\label{table3}
\end{table*}

\section{Exponentially improved asymptotic expansions}\label{section5}

We shall find it convenient to express our exponentially improved expansions in terms of the (scaled) Terminant function, which is defined by
\[
\widehat T_p \left( w \right) = \frac{{e^{\pi ip} w^{1 - p} e^{ - w} }}{{2\pi i}}\int_0^{ + \infty } {\frac{{t^{p - 1} e^{ - t} }}{w + t}dt} \; \text{ for } \; p>0 \; \text{ and } \; \left| \arg w \right| < \pi ,
\]
and by analytic continuation elsewhere. Olver \cite{Olver4} showed that when $p \sim \left|w\right|$ and $w \to \infty$, we have
\begin{equation}\label{eq36}
ie^{ - \pi ip} \widehat T_p \left( w \right) = \begin{cases} \mathcal{O}\left( {e^{ - w - \left| w \right|} } \right) & \; \text{ if } \; \left| {\arg w} \right| \le \pi \\ \mathcal{O}\left(1\right) & \; \text{ if } \; - 3\pi  < \arg w \le  - \pi. \end{cases}
\end{equation}
Concerning the smooth transition of the Stokes discontinuities, we will use the more precise asymptotic formulas
\begin{equation}\label{eq37}
\widehat T_p \left( w \right) = \frac{1}{2} + \frac{1}{2}\mathop{\text{erf}} \left( {c\left( \varphi  \right)\sqrt {\frac{1}{2}\left| w \right|} } \right) + \mathcal{O}\left( {\frac{{e^{ - \frac{1}{2}\left| w \right|c^2 \left( \varphi  \right)} }}{{\left| w \right|^{\frac{1}{2}} }}} \right)
\end{equation}
for $-\pi +\delta \leq \arg w \leq 3 \pi -\delta$, $0 < \delta  \le 2\pi$; and
\begin{equation}\label{eq38}
e^{ - 2\pi ip} \widehat T_p \left( w \right) =  - \frac{1}{2} + \frac{1}{2}\mathop{\text{erf}} \left( { - \overline {c\left( { - \varphi } \right)} \sqrt {\frac{1}{2}\left| w \right|} } \right) + \mathcal{O}\left( {\frac{{e^{ - \frac{1}{2}\left| w \right|\overline {c^2 \left( { - \varphi } \right)} } }}{{\left| w \right|^{\frac{1}{2}} }}} \right)
\end{equation}
for $- 3\pi  + \delta  \le \arg w \le \pi  - \delta$, $0 < \delta \le 2\pi$. Here $\varphi = \arg w$ and erf denotes the Error function. The quantity $c\left( \varphi  \right)$ is defined implicitly by the equation
\[
\frac{1}{2}c^2 \left( \varphi  \right) = 1 + i\left( {\varphi  - \pi } \right) - e^{i\left( {\varphi  - \pi } \right)},
\]
and corresponds to the branch of $c\left( \varphi  \right)$ which has the following expansion in the neighbourhood of $\varphi = \pi$:
\begin{equation}\label{eq39}
c\left( \varphi  \right) = \left( {\varphi  - \pi } \right) + \frac{i}{6}\left( {\varphi  - \pi } \right)^2  - \frac{1}{{36}}\left( {\varphi  - \pi } \right)^3  - \frac{i}{{270}}\left( {\varphi  - \pi } \right)^4  +  \cdots .
\end{equation}
For complete asymptotic expansions, see Olver \cite{Olver5}. We remark that Olver uses the different notation $F_p \left( w \right) = ie^{ - \pi ip} \widehat T_p \left( w \right)$ for the Terminant function and the other branch of the function $c\left( \varphi  \right)$. For further properties of the Terminant function, see, for example, Paris and Kaminski \cite[Chapter 6]{Paris3}.

\subsection{Proof of the exponentially improved expansions for $\mathbf{A}_{\nu}\left(\lambda \nu \right)$}

\subsubsection{Case (i): $\lambda>1$} The proof goes exactly the same way as the proof of Theorem 3 in the first paper of this series \cite{Nemes2}. One has to replace $R_N \left( {\nu ,\beta } \right)$, $R_{N,M} \left( {\nu ,\beta } \right)$ and $\tan \beta - \beta$ by $R_N \left( {\nu ,\sec \beta } \right)$, $R_{N,M} \left( {\nu ,\sec \beta } \right)$ and $\tan \beta - \beta + \pi$ in the corresponding formulas.

\subsubsection{Case (ii): $\lambda=1$}

First, we suppose that $\left|\arg \nu\right| < \frac{\pi}{2}$. Our starting point is the representation \eqref{eq10}, written in the form
\begin{equation}\label{eq40}
R_N \left( {\nu ,1 } \right) = \frac{{\left( { - 1} \right)^N }}{{2\pi \nu ^{2N + 1} }}\int_0^{ + \infty } {\frac{{t^{2N} e^{ - \pi t} }}{{1 - it/\nu }}iH_{it}^{\left( 1 \right)} \left( {it } \right)dt}  + \frac{{\left( { - 1} \right)^N }}{{2\pi \nu ^{2N + 1} }}\int_0^{ + \infty } {\frac{{t^{2N} e^{ - \pi t} }}{{1 + it/\nu }}iH_{it}^{\left( 1 \right)} \left( {it } \right)dt} .
\end{equation}
Let $0 \leq M <3N$ be a fixed integer such that $M \equiv 0 \mod 3$. We use \eqref{eq31} to expand the function $H_{it}^{\left( 1 \right)} \left( {it } \right)$ under the integrals in \eqref{eq40}, to obtain
\begin{gather}\label{eq41}
\begin{split}
R_N \left( {\nu ,1} \right) = \; & \frac{2}{{3\pi }}\sum\limits_{m = 0}^{M - 1} {d_{2m} \sin \left( {\frac{{\left( {2m + 1} \right)\pi }}{3}} \right)\frac{{\Gamma \left( {\frac{{2m + 1}}{3}} \right)}}{{\nu ^{\frac{{2m + 1}}{3}} }}\frac{{\left( { - 1} \right)^{N + m} \nu ^{\frac{{2m - 2}}{3} - 2N} }}{{2\pi }}\int_0^{ + \infty } {\frac{{t^{2N - \frac{{2m - 2}}{3} - 1} e^{ - \pi t} }}{{1 - it/\nu }}dt} }\\ 
& + \frac{2}{{3\pi }}\sum\limits_{m = 0}^{M - 1} {d_{2m} \sin \left( {\frac{{\left( {2m + 1} \right)\pi }}{3}} \right)\frac{{\Gamma \left( {\frac{{2m + 1}}{3}} \right)}}{{\nu ^{\frac{{2m + 1}}{3}} }}\frac{{\left( { - 1} \right)^{N + m } \nu ^{\frac{{2m - 2}}{3} - 2N} }}{{2\pi }}\int_0^{ + \infty } {\frac{{t^{2N - \frac{{2m - 2}}{3} - 1} e^{ - \pi t} }}{{1 + it/\nu }}dt} } \\
& + R_{N,M} \left( {\nu ,1} \right),
\end{split}
\end{gather}
with
\begin{equation}\label{eq42}
R_{N,M} \left( {\nu ,1} \right) = \frac{{\left( { - 1} \right)^N }}{{2\pi \nu ^{2N + 1} }}\int_0^{ + \infty } {\frac{{t^{2N} e^{ - \pi t} }}{{1 - it/\nu }}iR_M^{\left( H \right)} \left( {it} \right)dt}  + \frac{{\left( { - 1} \right)^N }}{{2\pi \nu ^{2N + 1} }}\int_0^{ + \infty } {\frac{{t^{2N} e^{ - \pi t} }}{{1 + it/\nu }}iR_M^{\left( H \right)} \left( {it} \right)dt} .
\end{equation}
The integrals in \eqref{eq41} can be identified in terms of the Terminant function since
\[
\frac{{\left( { - 1} \right)^{N + m} \nu ^{\frac{{2m - 2}}{3} - 2N} }}{{2\pi }}\int_0^{ + \infty } {\frac{{t^{2N - \frac{{2m - 2}}{3} - 1} e^{ - \pi t} }}{{1 - it/\nu }}dt}  = -ie^{\pi i\nu } e^{\frac{{2\left( {2m + 1} \right)\pi i}}{3}} \widehat T_{2N - \frac{{2m - 2}}{3}} \left( {\pi i\nu } \right)
\]
and
\[
\frac{{\left( { - 1} \right)^{N + m} \nu ^{\frac{{2m - 2}}{3} - 2N} }}{{2\pi }}\int_0^{ + \infty } {\frac{{t^{2N - \frac{{2m - 2}}{3} - 1} e^{ - \pi t} }}{{1 + it/\nu }}dt}  = -ie^{ - \pi i\nu } \widehat T_{2N - \frac{{2m - 2}}{3}} \left( { - \pi i\nu } \right).
\]
Hence, we have the following expansion
\begin{multline*}
R_N \left( {\nu ,1} \right) = -ie^{\pi i\nu } \frac{2}{{3\pi }}\sum\limits_{m = 0}^{M - 1} {d_{2m} e^{\frac{{2\left( {2m + 1} \right)\pi i}}{3}} \sin \left( {\frac{{\left( {2m + 1} \right)\pi }}{3}} \right)\frac{{\Gamma \left( {\frac{{2m + 1}}{3}} \right)}}{{\nu ^{\frac{{2m + 1}}{3}} }}\widehat T_{2N - \frac{{2m - 2}}{3}} \left( {\pi i\nu } \right)} 
\\ - ie^{ - \pi i\nu } \frac{2}{{3\pi }}\sum\limits_{m = 0}^{M - 1} {d_{2m} \sin \left( {\frac{{\left( {2m + 1} \right)\pi }}{3}} \right)\frac{{\Gamma \left( {\frac{{2m + 1}}{3}} \right)}}{{\nu ^{\frac{{2m + 1}}{3}} }}\widehat T_{2N - \frac{{2m - 2}}{3}} \left( { - \pi i\nu } \right)}  + R_{N,M} \left( {\nu ,1} \right).
\end{multline*}
Taking $\nu = r e^{i\theta}$, the representation \eqref{eq42} takes the form
\begin{equation}\label{eq43}
R_{N,M} \left( {\nu ,1} \right) = \frac{{\left( { - 1} \right)^N }}{{2\pi \left( {e^{i\theta } } \right)^{2N + 1} }}\int_0^{ + \infty } {\frac{{\tau ^{2N} e^{ - \pi r\tau } }}{{1 - i\tau e^{ - i\theta } }}iR_M^{\left( H \right)} \left( {ir\tau } \right)d\tau }  + \frac{{\left( { - 1} \right)^N }}{{2\pi \left( {e^{i\theta } } \right)^{2N + 1} }}\int_0^{ + \infty } {\frac{{\tau ^{2N} e^{ - \pi r\tau } }}{{1 + i\tau e^{ - i\theta } }}iR_M^{\left( H \right)} \left( {ir\tau } \right)d\tau } .
\end{equation}
In \cite[Appendix B]{Nemes} it was shown that
\[
\frac{{1 - \left( {s/r\tau } \right)^{\frac{4}{3}} }}{{1 - \left( {s/r\tau } \right)^2 }} = \frac{{1 - \left( {s/r} \right)^{\frac{4}{3}} }}{{1 - \left( {s/r} \right)^2 }} + \left( {\tau  - 1} \right)f\left( {r,\tau ,s} \right)
\]
for positive $r$, $\tau$ and $s$, with some $f\left(r,\tau ,s\right)$ satisfying $\left|f\left(r,\tau ,s\right)\right| \leq 2$. Using the integral formula \eqref{eq35}, $R_M^{\left( H \right)} \left( {ir\tau } \right)$ can be written as
\begin{align*}
R_M^{\left( H \right)} \left( {ir\tau } \right) = \; & \frac{1}{{\sqrt 3 \pi \left( {r\tau } \right)^{\frac{{2M + 1}}{3}} }}\int_0^{ + \infty } {s^{\frac{2M - 2}{3}} e^{ - 2\pi s} \frac{{1 - \left( {s/r\tau } \right)^{\frac{4}{3}} }}{{1 - \left( {s/r\tau } \right)^2 }}H_{is}^{\left( 1 \right)} \left( {is} \right)ds} \\
= \; & \frac{1}{{\sqrt 3 \pi \left( {r\tau } \right)^{\frac{{2M + 1}}{3}} }}\int_0^{ + \infty } {s^{\frac{{2M - 2}}{3}} e^{ - 2\pi s} \frac{{1 - \left( {s/r} \right)^{\frac{4}{3}} }}{{1 - \left( {s/r} \right)^2 }}H_{is}^{\left( 1 \right)} \left( {is} \right)ds} \\
& + \frac{{\tau  - 1}}{{\sqrt 3 \pi \left( {r\tau } \right)^{\frac{{2M + 1}}{3}} }}\int_0^{ + \infty } {s^{\frac{{2M - 2}}{3}} e^{ - 2\pi s} f\left( {r,\tau ,s} \right)H_{is}^{\left( 1 \right)} \left( {is} \right)ds} .
\end{align*}
Noting that
\[
0< \frac{{1 - \left( {s/r} \right)^{\frac{4}{3}} }}{{1 - \left( {s/r} \right)^2 }} < 1
\]
for any positive $r$ and $s$, substitution into \eqref{eq43} yields the upper bound
\begin{align*}
\left| {R_{N,M} \left( {\nu ,1} \right)} \right| \le \; &\frac{{\left| {d_{2M} } \right|\Gamma \left( {\frac{{2M + 1}}{3}} \right)}}{{\sqrt 3 \pi \left| \nu  \right|^{\frac{{2M + 1}}{3}} }}\left| {\frac{1}{{2\pi }}\int_0^{ + \infty } {\frac{{\tau ^{2N - \frac{{2M - 2}}{3} - 1} e^{ - \pi r\tau } }}{{1 - i\tau e^{ - i\theta } }}d\tau } } \right| \\ &+ \frac{{\left| {d_{2M} } \right|\Gamma \left( {\frac{{2M + 1}}{3}} \right)}}{{\sqrt 3 \pi ^2 \left| \nu  \right|^{\frac{{2M + 1}}{3}} }}\int_0^{ + \infty } {\tau ^{2N - \frac{{2M - 2}}{3} - 1} e^{ - \pi r\tau } \left| {\frac{{\tau  - 1}}{{\tau  + ie^{i\theta } }}} \right|d\tau } \\
& + \frac{{\left| {d_{2M} } \right|\Gamma \left( {\frac{{2M + 1}}{3}} \right)}}{{\sqrt 3 \pi \left| \nu  \right|^{\frac{{2M + 1}}{3}} }}\left| {\frac{1}{{2\pi }}\int_0^{ + \infty } {\frac{{\tau ^{2N - \frac{{2M - 2}}{3} - 1} e^{ - \pi r\tau } }}{{1 + i\tau e^{ - i\theta } }}d\tau } } \right| \\ & + \frac{{\left| {d_{2M} } \right|\Gamma \left( {\frac{{2M + 1}}{3}} \right)}}{{\sqrt 3 \pi ^2 \left| \nu  \right|^{\frac{{2M + 1}}{3}} }}\int_0^{ + \infty } {\tau ^{2N - \frac{{2M - 2}}{3} - 1} e^{ - \pi r\tau } \left| {\frac{{\tau  - 1}}{{\tau  - ie^{i\theta } }}} \right|d\tau } .
\end{align*}
As $\left| \left(\tau  - 1\right)/\left(\tau \pm i e^{i\theta}\right)  \right| \le 1$, we find that
\begin{align*}
\left| {R_{N,M} \left( {\nu ,1} \right)} \right| \le  \frac{{\left| {d_{2M} } \right|\Gamma \left( {\frac{{2M + 1}}{3}} \right)}}{{\sqrt 3 \pi \left| \nu  \right|^{\frac{{2M + 1}}{3}} }}\left| {e^{\pi i\nu } \widehat T_{2N - \frac{{2M - 2}}{3}} \left( {\pi i\nu } \right)} \right| & + \frac{{\left| {d_{2M} } \right|\Gamma \left( {\frac{{2M + 1}}{3}} \right)}}{{\sqrt 3 \pi \left| \nu  \right|^{\frac{{2M + 1}}{3}} }}\left| {e^{ - \pi i\nu } \widehat T_{2N - \frac{{2M - 2}}{3}} \left( { - \pi i\nu } \right)} \right|\\
& + \frac{{2\left| {d_{2M} } \right|\Gamma \left( {\frac{{2M + 1}}{3}} \right)\Gamma \left( {2N - \frac{{2M - 2}}{3}} \right)}}{{\sqrt 3 \pi ^2 \pi ^{2N - \frac{{2M - 2}}{3}} \left| \nu  \right|^{2N + 1} }}.
\end{align*}
By continuity, this bound holds in the closed sector $\left|\arg \nu\right| \le \frac{\pi}{2}$. Assume that $N = \frac{1}{2}\pi \left|\nu\right|+\rho$ where $\rho$ is bounded. Employing Stirling's formula, we find that
\[
\frac{{2\left| {d_{2M} } \right|\Gamma \left( {\frac{{2M + 1}}{3}} \right)\Gamma \left( {2N - \frac{{2M - 2}}{3}} \right)}}{{\sqrt 3 \pi ^2 \pi ^{2N - \frac{{2M - 2}}{3}} \left| \nu  \right|^{2N + 1} }} = \mathcal{O}_{M,\rho } \left( {\frac{{e^{ - \pi \left| \nu  \right|} }}{{\left| \nu  \right|^{\frac{1}{2}} }}\left| {d_{2M} } \right|\frac{{\Gamma \left( {\frac{{2M + 1}}{3}} \right)}}{{\left| \nu  \right|^{\frac{{2M + 1}}{3}} }}} \right)
\]
as $\nu \to \infty$. Olver's estimation \eqref{eq36} shows that
\[
\left| {e^{ \pm \pi i\nu } \widehat T_{2N - \frac{{2M - 2}}{3}} \left( { \pm \pi i\nu } \right)} \right| = \mathcal{O}_{M,\rho } \left( {e^{ - \pi \left| \nu  \right|} } \right)
\]
for large $\nu$. Therefore, we have that
\begin{equation}\label{eq44}
R_{N,M} \left( {\nu ,1} \right) = \mathcal{O}_{M,\rho } \left( {e^{ - \pi \left| \nu  \right|} \left| {d_{2M} } \right|\frac{{\Gamma \left( {\frac{{2M + 1}}{3}} \right)}}{{\left| \nu  \right|^{\frac{{2M + 1}}{3}} }}} \right)
\end{equation}
as $\nu \to \infty$ in the sector $\left|\arg \nu\right| \le \frac{\pi}{2}$.

Rotating the path of integration in \eqref{eq42} and applying the residue theorem yields
\begin{align*}
R_{N,M} \left( {\nu ,1} \right) & = ie^{\pi i\nu } R_M^{\left( H \right)} \left( \nu  \right) + \frac{{\left( { - 1} \right)^N }}{{2\pi \nu ^{2N + 1} }}\int_0^{ + \infty } {\frac{{t^{2N} e^{ - \pi t} }}{{1 - it/\nu }}iR_M^{\left( H \right)} \left( {it} \right)dt}  + \frac{{\left( { - 1} \right)^N }}{{2\pi \nu ^{2N + 1} }}\int_0^{ + \infty } {\frac{{t^{2N} e^{ - \pi t} }}{{1 + it/\nu }}iR_M^{\left( H \right)} \left( {it} \right)dt}\\
& = ie^{\pi i\nu } R_M^{\left( H \right)} \left( \nu  \right) - R_{N,M} \left( {\nu e^{ - \pi i} ,1} \right),
\end{align*}
when $\frac{\pi}{2} < \arg \nu < \frac{3\pi}{2}$. It follows that
\[
\left| {R_{N,M} \left( {\nu ,1} \right)} \right| \le e^{ - \pi \Im \left( \nu  \right)} \left| {R_M^{\left( H \right)} \left( \nu  \right)} \right| + \left| {R_{N,M} \left( {\nu e^{ - \pi i} ,1} \right)} \right|
\]
in the closed sector $\frac{\pi}{2} \leq \arg \nu \leq \frac{3\pi}{2}$, using continuity. It was proved in \cite{Nemes} that $R_M^{\left( H \right)} \left( \nu  \right) = \mathcal{O}_M \left( {\left| {d_{2M} } \right|\Gamma \left( {\frac{{2M + 1}}{3}} \right)\left| \nu  \right|^{ - \frac{{2M + 1}}{3}} } \right)$ as $\nu \to \infty$ in the closed sector $-\frac{\pi}{2} \leq \arg \nu \leq \frac{3\pi}{2}$, whence, by \eqref{eq44}, we deduce that
\begin{gather}\label{eq58}
\begin{split}
R_{N,M} \left( {\nu ,1} \right) & = \mathcal{O}_M \left( {e^{ - \pi \Im \left( \nu  \right)} \left| {d_{2M} } \right|\frac{{\Gamma \left( {\frac{{2M + 1}}{3}} \right)}}{{\left| \nu  \right|^{\frac{{2M + 1}}{3}} }}} \right) + \mathcal{O}_{M,\rho } \left( {e^{ - \pi \left| \nu  \right|} \left| {d_{2M} } \right|\frac{{\Gamma \left( {\frac{{2M + 1}}{3}} \right)}}{{\left| \nu  \right|^{\frac{{2M + 1}}{3}} }}} \right)
\\& = \mathcal{O}_{M,\rho } \left( {e^{ - \pi \Im \left( \nu  \right)} \left| {d_{2M} } \right|\frac{{\Gamma \left( {\frac{{2M + 1}}{3}} \right)}}{{\left| \nu  \right|^{\frac{{2M + 1}}{3}} }}} \right)
\end{split}
\end{gather}
as $\nu \to \infty$ in the sector $\frac{\pi}{2} \leq \arg \nu \leq \frac{3\pi}{2}$.

The reflection principle gives the relation
\[
R_{N,M} \left( {\nu ,1} \right) = \overline {R_{N,M} \left( {\bar \nu ,1} \right)}  =  - ie^{ - \pi i\nu } \overline {R_M^{\left( H \right)} \left( {\bar \nu } \right)}  - R_{N,M} \left( {\nu e^{\pi i} ,1} \right) =  - ie^{ - \pi i\nu } R_M^{\left( H \right)} \left( {\nu e^{\pi i} } \right) - R_{N,M} \left( {\nu e^{\pi i} ,1} \right),
\]
valid when $-\frac{3\pi}{2} < \arg \nu < -\frac{\pi}{2}$. Trivial estimation and a continuity argument show that
\[
\left| {R_{N,M} \left( {\nu ,1} \right)} \right| \le e^{\pi \Im \left( \nu  \right)} \left| {R_M^{\left( H \right)} \left( {\nu e^{\pi i} } \right)} \right| + \left| {R_{N,M} \left( {\nu e^{\pi i} ,1} \right)} \right|
\]
in the closed sector $-\frac{3\pi}{2} \leq \arg \nu \leq -\frac{\pi}{2}$. Since $R_M^{\left( H \right)} \left( \nu e^{\pi i} \right) = \mathcal{O}_M \left( {\left| {d_{2M} } \right|\Gamma \left( {\frac{{2M + 1}}{3}} \right)\left| \nu  \right|^{ - \frac{{2M + 1}}{3}} } \right)$ as $\nu \to \infty$ in the range $-\frac{3\pi}{2} \leq \arg \nu \leq -\frac{\pi}{2}$, by \eqref{eq44}, we find that
\begin{gather}\label{eq59}
\begin{split}
R_{N,M} \left( {\nu ,1} \right) & = \mathcal{O}_M \left( {e^{\pi \Im \left( \nu  \right)} \left| {d_{2M} } \right|\frac{{\Gamma \left( {\frac{{2M + 1}}{3}} \right)}}{{\left| \nu  \right|^{\frac{{2M + 1}}{3}} }}} \right) + \mathcal{O}_{M,\rho } \left( {e^{ - \pi \left| \nu  \right|} \left| {d_{2M} } \right|\frac{{\Gamma \left( {\frac{{2M + 1}}{3}} \right)}}{{\left| \nu  \right|^{\frac{{2M + 1}}{3}} }}} \right)
\\& = \mathcal{O}_{M,\rho } \left( {e^{\pi \Im \left( \nu  \right)} \left| {d_{2M} } \right|\frac{{\Gamma \left( {\frac{{2M + 1}}{3}} \right)}}{{\left| \nu  \right|^{\frac{{2M + 1}}{3}} }}} \right)
\end{split}
\end{gather}
as $\nu \to \infty$ with $-\frac{3\pi}{2} \leq \arg \nu \leq -\frac{\pi}{2}$.

If $M \equiv 1 \mod 3$ or $M \equiv 2 \mod 3$, we write the remainder $R_{N,M} \left( {\nu ,1} \right)$ in the form
\begin{align*}
R_{N,M} \left( {\nu ,1} \right) = & - ie^{\pi i\nu } \frac{2}{{3\pi }}d_{2M + 2} e^{\frac{\pi }{3}i} \frac{{\sqrt 3 }}{2}\frac{{\Gamma \left( {\frac{{2M + 3}}{3}} \right)}}{{\nu ^{\frac{{2M + 3}}{3}} }}\widehat T_{2N - \frac{{2M}}{3}} \left( {\pi i\nu } \right)\\
& + ie^{ - \pi i\nu } \frac{2}{{3\pi }}d_{2M + 2} \frac{{\sqrt 3 }}{2}\frac{{\Gamma \left( {\frac{{2M + 3}}{3}} \right)}}{{\nu ^{\frac{{2M + 3}}{3}} }}\widehat T_{2N - \frac{{2M}}{3}} \left( { - \pi i\nu } \right) + R_{N,M + 2} \left( {\nu ,1} \right)
\end{align*}
and
\begin{align*}
R_{N,M} \left( {\nu ,1} \right) = & - ie^{\pi i\nu } \frac{2}{{3\pi }}d_{2M} e^{\frac{\pi }{3}i} \frac{{\sqrt 3 }}{2}\frac{{\Gamma \left( {\frac{{2M + 1}}{3}} \right)}}{{\nu ^{\frac{{2M + 1}}{3}} }}\widehat T_{2N - \frac{{2M - 2}}{3}} \left( {\pi i\nu } \right)\\
& + ie^{ - \pi i\nu } \frac{2}{{3\pi }}d_{2M} \frac{{\sqrt 3 }}{2}\frac{{\Gamma \left( {\frac{{2M + 1}}{3}} \right)}}{{\nu ^{\frac{{2M + 1}}{3}} }}\widehat T_{2N - \frac{{2M - 2}}{3}} \left( { - \pi i\nu } \right) + R_{N,M + 1} \left( {\nu ,1} \right),
\end{align*}
respectively. Applying the connection formula $\widehat T_p \left( w \right) = e^{2\pi ip} \left( {\widehat T_p \left( {we^{2\pi i} } \right) - 1} \right)$ together with Olver's result \eqref{eq36} and the bounds \eqref{eq44}, \eqref{eq58}, \eqref{eq59} we have established, the estimates for the cases $M \equiv 1 \mod 3$ and $M \equiv 2 \mod 3$ follow.

\subsection{Stokes phenomenon and Berry's transition}

\subsubsection{Case (i): $\lambda>1$} As usual, let $\lambda = \sec \beta$ with some $0 <\beta <\frac{\pi}{2}$. We study the Stokes phenomenon related to the asymptotic expansion of $\mathbf{A}_{\nu} \left(\nu \sec \beta\right)$ occurring when $\arg \nu$ passes through the values $\pm \frac{\pi}{2}$. In the range $\left|\arg \nu\right|<\frac{\pi}{2}$, the asymptotic expansion
\begin{equation}\label{eq47}
\mathbf{A}_{\nu} \left( {\nu \sec \beta } \right) \sim  \frac{1}{\pi }\sum\limits_{n = 0}^\infty  {\frac{{\left( {2n} \right)!a_n \left(\sec \beta\right)}}{{\nu ^{2n + 1} }}} 
\end{equation}
holds as $\nu \to \infty$. From \eqref{eq50} we have
\[
\mathbf{A}_\nu  \left( {\nu \sec \beta } \right) = i e^{\pi i\nu } H_\nu ^{\left( 1 \right)} \left( {\nu \sec \beta } \right) - \mathbf{A}_{-\nu} \left( {\nu e^{ - \pi i} \sec \beta } \right)
\]
when $\frac{\pi}{2} < \arg \nu < \frac{3\pi}{2}$, and
\[
\mathbf{A}_\nu  \left( {\nu \sec \beta } \right) =  - ie^{ - \pi i\nu } H_\nu ^{\left( 2 \right)} \left( {\nu \sec \beta } \right) - \mathbf{A}_{-\nu} \left( {\nu e^{\pi i} \sec \beta } \right)
\]
for $-\frac{3\pi}{2} < \arg \nu < -\frac{\pi}{2}$. For the right-hand sides, we can apply the asymptotic expansions of the Hankel functions and the Anger--Weber function to deduce that
\begin{equation}\label{eq54}
\mathbf{A}_{\nu } \left( {\nu \sec \beta } \right) \sim i\frac{{e^{i\nu \left( {\tan \beta  - \beta +\pi} \right) - \frac{\pi }{4}i} }}{{\left( {\frac{1}{2}\nu \pi \tan \beta } \right)^{\frac{1}{2}} }}\sum\limits_{m = 0}^\infty  {\left( { - 1} \right)^m \frac{{U_m \left( {i\cot \beta } \right)}}{{\nu ^m }}}  + \frac{1}{\pi}\sum\limits_{n = 0}^\infty  {\frac{{\left( {2n} \right)!a_n \left( {\sec \beta } \right)}}{{\nu ^{2n + 1} }}} 
\end{equation}
as $\nu \to \infty$ in the sector $\frac{\pi}{2} < \arg \nu < \frac{3\pi}{2}$, and
\begin{equation}\label{eq55}
\mathbf{A}_{\nu } \left( {\nu \sec \beta } \right) \sim  - i\frac{{e^{ - i\nu \left( {\tan \beta  - \beta +\pi} \right) + \frac{\pi }{4}i} }}{{\left( {\frac{1}{2}\nu \pi \tan \beta } \right)^{\frac{1}{2}} }}\sum\limits_{m = 0}^\infty  {\frac{{U_m \left( {i\cot \beta } \right)}}{{\nu ^m }}}  + \frac{1}{\pi}\sum\limits_{n = 0}^\infty  {\frac{{\left( {2n} \right)!a_n \left( {\sec \beta } \right)}}{{\nu ^{2n + 1} }}} 
\end{equation}
as $\nu \to \infty$ in the sector $-\frac{3\pi}{2} < \arg \nu < -\frac{\pi}{2}$. Therefore, as the line $\arg \nu = \frac{\pi}{2}$ is crossed, the additional series
\begin{equation}\label{eq48}
i\frac{{e^{i\nu \left( {\tan \beta  - \beta +\pi} \right) - \frac{\pi }{4}i} }}{{\left( {\frac{1}{2}\nu \pi \tan \beta } \right)^{\frac{1}{2}} }}\sum\limits_{m = 0}^\infty  {\left( { - 1} \right)^m \frac{{U_m \left( {i\cot \beta } \right)}}{\nu ^m}}
\end{equation}
appears in the asymptotic expansion of $\mathbf{A}_{\nu} \left( {\nu \sec \beta } \right)$ beside the original one \eqref{eq47}. Similarly, as we pass through the line $\arg \nu = -\frac{\pi}{2}$, the series
\begin{equation}\label{eq49}
- i\frac{{e^{ - i\nu \left( {\tan \beta  - \beta +\pi} \right) + \frac{\pi }{4}i} }}{{\left( {\frac{1}{2}\nu \pi \tan \beta } \right)^{\frac{1}{2}} }}\sum\limits_{m = 0}^\infty  {\frac{{U_m \left( {i\cot \beta } \right)}}{\nu ^m }}
\end{equation}
appears in the asymptotic expansion of $\mathbf{A}_{\nu} \left( {\nu \sec \beta } \right)$ beside the original series \eqref{eq47}. We have encountered a Stokes phenomenon with Stokes lines $\arg \nu = \pm\frac{\pi}{2}$.

In his important paper \cite{Berry2}, Berry gave a new interpretation of the Stokes phenomenon; he proved that assuming optimal truncation, the transition between compound asymptotic expansions is of Error function type, thus yielding a smooth and rapid transition as a Stokes line is crossed.

Using the exponentially improved expansion given in Theorem \ref{thm2}, we show that the asymptotic expansion of $\mathbf{A}_{\nu} \left( {\nu \sec \beta } \right)$ exhibits the Berry transition between the two asymptotic series across the Stokes lines $\arg \nu = \pm\frac{\pi}{2}$. More precisely, we shall find that the first few terms of the series in \eqref{eq48} and \eqref{eq49} ``emerge" in a rapid and smooth way as $\arg \nu$ passes through $\frac{\pi}{2}$ and $-\frac{\pi}{2}$, respectively.

From Theorem \ref{thm2}, we conclude that if $N \approx \frac{1}{2}\left| \nu  \right|\left( {\tan \beta  - \beta +\pi} \right)$, then for large $\nu$, $ \left|\arg \nu\right| < \pi$, we have
\begin{align*}
\mathbf{A}_{\nu} \left( {\nu \sec \beta } \right) \approx \; & \frac{1}{\pi }\sum\limits_{n = 0}^{N - 1} {\frac{{\left( {2n} \right)!a_n \left( {\sec \beta } \right)}}{\nu ^{2n + 1}}} \\ & + i\frac{{e^{i\nu \left( {\tan \beta  - \beta +\pi} \right) - \frac{\pi }{4}i} }}{{\left( {\frac{1}{2}\nu \pi \tan \beta } \right)^{\frac{1}{2}} }}\sum\limits_{m = 0} {\left( { - 1} \right)^m \frac{{U_m \left( {i\cot \beta } \right)}}{{\nu ^m }}\widehat T_{2N - m + \frac{1}{2}} \left( {i\nu \left( {\tan \beta  - \beta +\pi} \right)} \right)} 
\\ & - i\frac{{e^{ - i\nu \left( {\tan \beta  - \beta +\pi} \right) + \frac{\pi }{4}i} }}{{\left( {\frac{1}{2}\nu \pi \tan \beta } \right)^{\frac{1}{2}} }}\sum\limits_{m = 0} {\frac{{U_m \left( {i\cot \beta } \right)}}{{\nu ^m }}\widehat T_{2N - m + \frac{1}{2}} \left( { - i\nu \left( {\tan \beta  - \beta +\pi} \right)} \right)} ,
\end{align*}
where $\sum\nolimits_{m = 0}$ means that the sum is restricted to the leading terms of the series.

In the upper half-plane, the terms involving $\widehat T_{2N - m + \frac{1}{2}} \left( { - i\nu \left( {\tan \beta  - \beta +\pi} \right)} \right)$ are exponentially small, the dominant contribution comes from the terms involving $\widehat T_{2N - m + \frac{1}{2}} \left( {i\nu \left( {\tan \beta  - \beta +\pi} \right)} \right)$. Under the above assumption on $N$, from \eqref{eq37} and \eqref{eq39}, the Terminant functions have the asymptotic behaviour
\[
\widehat T_{2N - m + \frac{1}{2}} \left( {i\nu \left( {\tan \beta  - \beta +\pi} \right)} \right) \sim \frac{1}{2} + \frac{1}{2} \mathop{\text{erf}} \left( {\left( {\theta  - \frac{\pi }{2}} \right)\sqrt {\frac{1}{2}\left| \nu  \right|\left( {\tan \beta  - \beta +\pi} \right)} } \right)
\]
provided that $\arg \nu = \theta$ is close to $\frac{\pi}{2}$, $\nu$ is large and $m$ is small in comparison with $N$. Therefore, when $\theta  < \frac{\pi}{2}$, the Terminant functions are exponentially small; for $\theta  = \frac{\pi }{2}$, they are asymptotically $\frac{1}{2}$ up to an exponentially small error; and when $\theta  >  \frac{\pi}{2}$, the Terminant functions are asymptotic to $1$ with an exponentially small error. Thus, the transition across the Stokes line $\arg \nu = \frac{\pi}{2}$ is effected rapidly and smoothly. Similarly, in the lower half-plane, the dominant contribution is controlled by the terms involving $\widehat T_{2N - m + \frac{1}{2}} \left( { - i\nu \left( {\tan \beta  - \beta +\pi} \right)} \right)$. From \eqref{eq38} and \eqref{eq39}, we have
\[
\widehat T_{2N - m + \frac{1}{2}} \left( { - i\nu \left( {\tan \beta  - \beta +\pi} \right)} \right) \sim \frac{1}{2} - \frac{1}{2} \mathop{\text{erf}} \left( {\left( {\theta  + \frac{\pi }{2}} \right)\sqrt {\frac{1}{2}\left| \nu  \right|\left( {\tan \beta  - \beta +\pi} \right)} } \right)
\]
under the assumptions that $\arg \nu = \theta$ is close to $-\frac{\pi}{2}$, $\nu$ is large and $m$ is small in comparison with $N \approx \frac{1}{2}\left| \nu  \right|\left( {\tan \beta  - \beta +\pi} \right)$. Thus, when $\theta  >  - \frac{\pi}{2}$, the Terminant functions are exponentially small; for $\theta  =  -\frac{\pi}{2}$, they are asymptotic to $\frac{1}{2}$ with an exponentially small error; and when $\theta < - \frac{\pi}{2}$, the Terminant functions are asymptotically $1$ up to an exponentially small error. Therefore, the transition through the Stokes line $\arg \nu = -\frac{\pi}{2}$ is carried out rapidly and smoothly.

We remark that from the expansions \eqref{eq54} and \eqref{eq55}, it follows that \eqref{eq47} is an asymptotic expansion of $\mathbf{A}_{\nu} \left( {\nu \sec \beta } \right)$ in the wider sector $\left|\arg \nu\right| \leq \pi -\delta < \pi$, with any fixed $0 < \delta \leq \pi$.

\subsubsection{Case (ii): $\lambda=1$} The analysis of the Stokes phenomenon for the asymptotic expansion of
$\mathbf{A}_{\nu} \left( {\nu} \right)$ is similar to the case $\lambda > 1$. In the range $\left|\arg \nu\right| <\frac{\pi}{2}$, the asymptotic expansion
\begin{equation}\label{eq51}
\mathbf{A}_{\nu} \left(\nu\right) \sim \frac{1}{\pi }\sum\limits_{n = 0}^\infty  {\frac{{\left( {2n} \right)!a_n \left(1\right)}}{{\nu ^{2n + 1} }}} 
\end{equation}
holds as $\nu \to \infty$. Employing the continuation formulas stated in Section \ref{section1}, we find that
\[
\mathbf{A}_\nu  \left( {\nu } \right) = i e^{\pi i\nu } H_\nu ^{\left( 1 \right)} \left( {\nu } \right) - \mathbf{A}_{-\nu} \left( {\nu e^{ - \pi i}} \right)
\]
and
\[
\mathbf{A}_\nu  \left( {\nu } \right) =  - ie^{ - \pi i\nu } H_\nu ^{\left( 2 \right)} \left( {\nu } \right) - \mathbf{A}_{-\nu} \left( {\nu e^{\pi i} } \right).
\]
For the right-hand sides, we can apply the asymptotic expansions of the Hankel functions and the Anger--Weber function to deduce that
\begin{equation}\label{eq56}
\mathbf{A}_{\nu} \left(\nu\right) \sim - ie^{\pi i \nu} \frac{2}{{3\pi }}\sum\limits_{m = 0}^\infty  {d_{2m} e^{\frac{{2\left( {2m + 1} \right)\pi i}}{3}} \sin \left( {\frac{{\left( {2m + 1} \right)\pi }}{3}} \right)\frac{{\Gamma \left( {\frac{{2m + 1}}{3}} \right)}}{{\nu ^{\frac{{2m + 1}}{3}} }}}  + \frac{1}{\pi }\sum\limits_{n = 0}^\infty  {\frac{{\left( {2n} \right)!a_n \left( 1 \right)}}{{\nu ^{2n + 1} }}} 
\end{equation}
as $\nu \to \infty$ in the sector $\frac{\pi}{2}<\arg \nu < \frac{3\pi}{2}$, and
\begin{equation}\label{eq57}
\mathbf{A}_{\nu} \left(\nu\right) \sim ie^{ - \pi i \nu} \frac{2}{{3\pi }}\sum\limits_{m = 0}^\infty  {d_{2m} e^{ - \frac{{2\left( {2m + 1} \right)\pi i}}{3}} \sin \left( {\frac{{\left( {2m + 1} \right)\pi }}{3}} \right)\frac{{\Gamma \left( {\frac{{2m + 1}}{3}} \right)}}{{\nu ^{\frac{{2m + 1}}{3}} }}}  + \frac{1}{\pi }\sum\limits_{n = 0}^\infty  {\frac{{\left( {2n} \right)!a_n \left( 1 \right)}}{{\nu ^{2n + 1} }}} 
\end{equation}
as $\nu \to \infty$ in the sector $-\frac{3\pi}{2}<\arg \nu < -\frac{\pi}{2}$. Therefore, as the line $\arg \nu = \frac{\pi}{2}$ is crossed, the additional series
\begin{equation}\label{eq52}
- ie^{\pi i \nu} \frac{2}{{3\pi }}\sum\limits_{m = 0}^\infty  {d_{2m} e^{\frac{{2\left( {2m + 1} \right)\pi i}}{3}} \sin \left( {\frac{{\left( {2m + 1} \right)\pi }}{3}} \right)\frac{{\Gamma \left( {\frac{{2m + 1}}{3}} \right)}}{{\nu ^{\frac{{2m + 1}}{3}} }}} 
\end{equation}
appears in the asymptotic expansion of $\mathbf{A}_{\nu} \left(\nu\right)$ beside the original one \eqref{eq51}. Similarly, as we pass through the line $\arg \nu = -\frac{\pi}{2}$, the series
\begin{equation}\label{eq53}
ie^{ - \pi i \nu} \frac{2}{{3\pi }}\sum\limits_{m = 0}^\infty  {d_{2m} e^{ - \frac{{2\left( {2m + 1} \right)\pi i}}{3}} \sin \left( {\frac{{\left( {2m + 1} \right)\pi }}{3}} \right)\frac{{\Gamma \left( {\frac{{2m + 1}}{3}} \right)}}{{\nu ^{\frac{{2m + 1}}{3}} }}}
\end{equation}
appears in the asymptotic expansion of $\mathbf{A}_{\nu} \left(\nu\right)$ beside the original series \eqref{eq51}. We have encountered a Stokes phenomenon with Stokes lines $\arg \nu = \pm\frac{\pi}{2}$. With the aid of the exponentially improved expansion given in Theorem \ref{thm3}, we shall find that the asymptotic series of $\mathbf{A}_{\nu} \left( \nu  \right)$ shows the Berry transition property: the two series in \eqref{eq52} and \eqref{eq53}  "emerge" in a rapid and smooth way as the Stokes lines $\arg \nu = \frac{\pi}{2}$ and $\arg \nu = -\frac{\pi}{2}$ are crossed.

From Theorem \ref{thm3}, we infer that if $N \approx \frac{1}{2}\pi \left|\nu\right|$, then for large $\nu$, $\left|\arg \nu\right|<\pi$, we have
\begin{multline*}
\mathbf{A}_\nu  \left( \nu  \right) \approx \frac{1}{\pi }\sum\limits_{n = 0}^{N - 1} {\frac{{\left( {2n} \right)!a_n \left( 1 \right)}}{{\nu ^{2n + 1} }}}  - ie^{\pi i\nu } \frac{2}{{3\pi }}\sum\limits_{m = 0} {d_{2m} e^{\frac{{2\left( {2m + 1} \right)\pi i}}{3}} \sin \left( {\frac{{\left( {2m + 1} \right)\pi }}{3}} \right)\frac{{\Gamma \left( {\frac{{2m + 1}}{3}} \right)}}{{\nu ^{\frac{{2m + 1}}{3}} }}\widehat T_{2N - \frac{{2m - 2}}{3}} \left( {\pi i\nu } \right)} \\ - ie^{ - \pi i\nu } \frac{2}{{3\pi }}\sum\limits_{m = 0} {d_{2m} e^{ - \frac{{2\left( {2m + 1} \right)\pi i}}{3}} \sin \left( {\frac{{\left( {2m + 1} \right)\pi }}{3}} \right)\frac{{\Gamma \left( {\frac{{2m + 1}}{3}} \right)}}{{\nu ^{\frac{{2m + 1}}{3}} }}e^{\frac{{2\left( {2m - 2} \right)\pi i}}{3}} \widehat T_{2N - \frac{{2m - 2}}{3}} \left( { - \pi i\nu } \right)} ,
\end{multline*}
where, as before, $\sum\nolimits_{m = 0}$ means that the sum is restricted to the leading terms of the series.

In the upper half-plane, the main contribution comes from the terms involving $\widehat T_{2N - \frac{{2m - 2}}{3}} \left( {\pi i\nu } \right)$. Under the above assumption on $N$, from \eqref{eq37} and \eqref{eq39}, the Terminant functions have the asymptotic behaviour
\[
\widehat T_{2N - \frac{{2m - 2}}{3}} \left( {\pi i\nu } \right) \sim \frac{1}{2} + \frac{1}{2}\mathop{\text{erf}} \left( {\left( {\theta  - \frac{\pi }{2}} \right)\sqrt {\frac{1}{2}\pi \left| \nu  \right|} } \right),
\]
provided that $\arg \nu = \theta$ is close to $\frac{\pi}{2}$, $\nu$ is large and $m$ is small in comparison with $N$. Therefore, when $\theta  < \frac{\pi}{2}$, the Terminant functions are exponentially small; for $\theta  = \frac{\pi }{2}$, they are asymptotically $\frac{1}{2}$ up to an exponentially small error; and when $\theta  >  \frac{\pi}{2}$, the Terminant functions are asymptotic to $1$ with an exponentially small error. Thus, the transition across the Stokes line $\arg \nu = \frac{\pi}{2}$ is effected rapidly and smoothly. Similarly, in the lower half-plane, the dominant contribution is controlled by the terms containing $\widehat T_{2N - \frac{{2m - 2}}{3}} \left( { - \pi i\nu } \right)$. From \eqref{eq38} and \eqref{eq39}, we have
\[
e^{\frac{{2\left( {2m - 2} \right)\pi i}}{3}} \widehat T_{2N - \frac{{2m - 2}}{3}} \left( { - \pi i\nu } \right) \sim  - \frac{1}{2} + \frac{1}{2}\mathop{\text{erf}}\left( {\left( {\theta  + \frac{\pi }{2}} \right)\sqrt {\frac{1}{2}\pi \left| \nu  \right|} } \right),
\]
under the assumptions that $\arg \nu = \theta$ is close to $-\frac{\pi}{2}$, $\nu$ is large and $m$ is small in comparison with $N \approx \frac{1}{2}\pi \left|\nu\right|$. Thus, when $\theta  >  - \frac{\pi}{2}$, the normalized Terminant functions are exponentially small; for $\theta  =  -\frac{\pi}{2}$, they are asymptotic to $-\frac{1}{2}$ with an exponentially small error; and when $\theta < - \frac{\pi}{2}$, the normalized Terminant functions are asymptotically $-1$ up to an exponentially small error. Therefore, the transition through the Stokes line $\arg \nu = -\frac{\pi}{2}$ is carried out rapidly and smoothly.

We note that from the expansions \eqref{eq56} and \eqref{eq57}, it follows that \eqref{eq51} is an asymptotic series of $\mathbf{A}_{\nu} \left( \nu \right)$ in the wider range $\left|\arg \nu\right| \leq \pi -\delta < \pi$, with any fixed $0 < \delta \leq \pi$.

\section{Discussion}\label{section6}

In this paper, we have discussed in detail the large order and argument asymptotics of the Anger--Weber function $\mathbf{A}_{\nu}\left(\lambda \nu\right)$ when $\lambda > 0$, using Howls' method. The resurgence properties and the exponentially improved versions of the large $\nu$ asymptotics of the associated Anger function $\mathbf{J}_\nu \left(\lambda \nu\right)$ and Weber function $\mathbf{E}_\nu \left(\lambda \nu\right)$ can be obtained from the relations
\[
\mathbf{J}_\nu  \left( {\lambda \nu } \right) = J_\nu  \left( {\lambda \nu } \right) + \sin \left( {\pi \nu } \right)\mathbf{A}_\nu  \left( {\lambda \nu } \right),
\]
\[
\mathbf{E}_\nu  \left( {\lambda \nu } \right) =  - Y_\nu  \left( {\lambda \nu } \right) - \cos \left( {\pi \nu } \right)\mathbf{A}_\nu  \left( {\lambda \nu } \right) - \mathbf{A}_{ - \nu } \left( {\lambda \nu } \right),
\]
our previous results on the Bessel functions \cite{Nemes} and the results of the present series of papers on the Anger--Weber function. Note that the resulting resurgence formulas have different forms according to whether $\lambda=1$ or $\lambda>1$. From these new representations, error bounds for the asymptotic expansions of the Anger and Weber functions can be derived which, in the case $\lambda=1$, may be compared with the ones given earlier by Olver \cite{Olver2}.

\section*{Acknowledgement} I would like to thank the anonymous referee for his/her constructive and helpful comments and suggestions on the manuscript.

\end{document}